\documentclass[review]{elsarticle}

\usepackage{lipsum}
\usepackage{amsfonts}
\usepackage{graphicx}
\usepackage{epstopdf}
\usepackage{algorithmic}
\usepackage{changes}
\usepackage[caption=false]{subfig}
\biboptions{numbers,sort&compress}

\title{Quantifying dependencies for sensitivity analysis with multivariate input sample data}
\author[cwi,cor]{A.W. ~Eggels}
\ead{a.w.eggels@cwi.nl}
\author[cwi,uva]{D.T. ~Crommelin}
\ead{daan.crommelin@cwi.nl}

\cortext[cor1]{Corresponding author.}
\address[cwi]{Centrum voor Wiskunde \& Informatica, P.O. Box 94079, 1090 GB Amsterdam, the Netherlands}
\address[uva]{Korteweg-de Vries institute for mathematics, University of Amsterdam, P.O. Box 94248, 1090 GE Amsterdam, the Netherlands}

\begin{document}
\begin{frontmatter}

\begin{abstract}
We present a novel method for quantifying dependencies in multivariate datasets, based on estimating the R\'{e}nyi entropy by minimum spanning trees (MSTs).
The length of the MSTs can be used to order pairs of variables from strongly to weakly dependent, making it a useful tool for sensitivity analysis with dependent input variables. It is well-suited for cases where the input distribution is unknown and only a sample of the inputs is available.
We introduce an estimator to quantify dependency based on the MST length, and investigate its properties with several numerical examples.
To reduce the computational cost of constructing the exact MST for large datasets, we explore methods to compute approximations to the exact MST, and find the multilevel approach introduced recently by Zhong et al. (2015) to be the most accurate. We apply our proposed method to an artificial testcase based on the Ishigami function, as well as to a real-world testcase involving sediment transport in the North Sea. The results are consistent with prior knowledge and heuristic understanding, as well as with variance-based analysis using Sobol indices in the case where these indices can be computed.
\end{abstract}

\end{frontmatter}
\section{Introduction}
Sensitivity analysis (SA) is a core topic in the field of Uncertainty Quantification (UQ) \cite{Mai10,2Sul15b}, in which the uncertainties in simulation models of complex systems are explored. SA addresses the question how uncertainty in the output of a model can be allocated to different uncertain model inputs \cite{2Cre09,Hel06,Hom96,2Ioo15,2Krz06,2Liu06,2Oak04,Sal02b,2Sal10,Sal08,Sob93,Sud08,Tru06,Wag95,Yam05}. Being able to answer this question allows one to focus on model inputs that induce the largest uncertainties in the output,  bringing benefits such as more efficient numerical exploration of model uncertainties, or guidance on how to reduce model output uncertainty most effectively.

A variety of methods for SA exists, both local and global methods, see e.g. the chapter on SA in \cite{3Gha17}, as well as the references mentioned above, for an overview. Many of these methods are intended for cases where the inputs are mutually independent, however there are situations where the assumption of independent inputs is unrealistic.

Input dependencies give several complications when using well-established methods. In that case, most methods either fail or make restrictive assumptions on the input, the output function, the number of available simulations, or combinations hereof. Examples of restrictions include assuming linear dependencies, normally distributed input variables, or being able to perform enough simulations to allow for Monte-Carlo methods \cite{2Cha15,2Cha12,2Dav09,2Dou13,2Kuc12,2Mar12,2Xu08}. The most well-known method is using Sobol's indices \cite{Sob93,3Sob01}. In case of dependent inputs, the functions in the Sobol decomposition are not orthogonal. Due to this, the indices may become negative, and therefore the interpretation of the indices becomes less clear \cite{Nav14}. Also, the contribution due to dependency can cancel the contribution due to the variable itself.

The setting under consideration in this paper is one in which the probability distribution of the multivariate input is unknown (both the joint distribution and the marginals), and only a sample (dataset) of the input is available. Thus, not only the dependency structure  (i.e., between which variables dependencies exist) is unknown, but also the shape of these dependencies. Our aim in this study is to develop a methodology to detect and quantify the dependencies and their strength from the input data, without making assumptions on the shape.

\cite{3Fil17} developed a Bayesian nonparametric procedure that leads to an analytic quantification for dependence versus independence. However, this procedure does not quantify explicitly the strength of the dependency. Furthermore, the inference resulting from the P\'{o}lya trees used in the procedure is known to depend strongly on the choice of the partitioning of the data points, although a partial optimization of the partitioning is proposed as well. This optimization increases the evidence in favour of dependence, but this is not a problem as long as the order of pairs of variables is considered in contrast to the numerical value of the evidence.

We build our methodology on the concept of entropy. Entropy is a notion that is used in various fields, e.g. statistical physics, information theory and mathematics, and that can represent dependencies between variables that go beyond linear relationships (correlation).
Its usefulness for SA was discussed before in \cite{2Aud08,2Liu06}. However, although entropy has a clear mathematical foundation, the estimation of entropy from data sampled from a continuous distribution is not straightforward. Both the use of binning (as in \cite{2Aud08}) and the Kullback-Leibler entropy (as in \cite{2Liu06}) are difficult, as the results with binning can be sensitive to the choice of bins, while the method to compute the Kullback-Leibler entropy can be sensitive to the parameters of the kernel-density estimation that it uses.

An appealing and elegant method for estimating entropies is due to Hero et al.\ \cite{2Her03,2Her02,2Her98}, who proposed to compute minimum spanning trees (MSTs) of the data and to use the length of the MST in an estimator for R\'{e}nyi entropy.
In this study, we employ this method for the purpose of SA, by using the length of the MST as a measure of dependence between variables. In principle, MSTs can be computed exactly, however the computation becomes expensive in case of large datasets. Therefore we test several approximate methods to compute the MST length with reduced computational cost.

In Section \ref{sec:setup} we briefly review some theory regarding entropy and estimators of entropy. Afterwards, in Section \ref{sec:methods} the approximation methods are discussed. Section \ref{sec:results} consists of validation of the proposed estimator and determining the consistency and robustness of the approximations to it. Section \ref{sec:5} links dependency quantification to sensitivity analysis, of which Section \ref{sec:tests} describes one example and one application. Section \ref{sec:discussion} concludes.

\section{Dependence and entropy}\label{sec:setup}
In this section we summarize a few basic definitions and concepts regarding dependence and entropy, as well as the connection between the two. This paves the way for defining our proposed estimator for dependency between two variables.

\subsection{Entropy}
Entropy is a concept which quantifies uncertainty or randomness in a system. It has its origins in classical thermodynamics, but generalizations to other fields, including mathematics, can be made.
Shannon's entropy (information theory) is defined as follows \cite{2Cov06}:
\begin{equation}\label{eq:Shannon}
H(X) = - \sum_{i} p_i \log(p_i),
\end{equation}
where $i$ represents the possible states of the discrete random variable $X$ and $p_i$ their probability of occurring. Its continuous analogue (differential entropy) is given by \cite{2Cov06}
\begin{equation}\label{eq:Shannoncontinuous}
H(X) = - \int_\Omega p(x)\log(p(x))dx,
\end{equation}
where $\Omega$ is the domain of the continuous random variable $X$. R\'{e}nyi has extended Shannon entropy to
\begin{equation}\label{eq:Renyi}
H_\alpha(X) = \frac{1}{1-\alpha} \log\left(\int_\Omega \left(p(x)\right)^\alpha dx\right),
\end{equation}
for $\alpha\in(0,\infty)$, see \cite{2Cov06}. In the limit $\alpha\rightarrow 1$, the R\'{e}nyi entropy (\ref{eq:Renyi}), sometimes referred to as  $\alpha$-entropy, reduces to the differential entropy (\ref{eq:Shannoncontinuous}).  For large values of $\alpha$, the events with high probability density determine the value of the entropy, whereas different values of the probability density are weighted more equally when $\alpha$ is close to $0$.

\subsection{Entropy as measure of dependence}
It is well-known that the dependence between two random variables $Y$ and $Z$ can be characterized by the mutual information (MI), i.e. the difference between the Shannon entropy of the joint distribution and the sum of the Shannon entropies of the marginals,
\begin{equation}\label{eq:MI}
I(Y,Z) = H(Y) + H(Z) - H(Y,Z) \, ,
\end{equation}
see \cite{2Cov06}. $I(Y,Z)$ attains it maximum if $Y$ and $Z$ are completely dependent, in that case $H(Y,Z)=H(Y)=H(Z)$ and $I(X,Y)=H(Y)$. Thus, $I(Y,Z)$ quantifies the strength of dependencies: the higher $I(Y,Z)$, the stronger the dependency.

In a similar manner, dependence can be quantified in terms of the R\'enyi entropy of the joint distribution $(Y,Z)$. The difference $H_\alpha(Y)+H_\alpha(Z) - H_\alpha(Y,Z)$ can be negative (unlike mutual information), however if we eliminate the effect of the marginal distributions we can use $H_\alpha(Y,Z)$ to quantify dependency. The lower $H_\alpha(Y,Z)$, the stronger the dependence. The marginals need to be eliminated because the R\'{e}nyi entropy is not scale-invariant. Therefore, the value of $H_\alpha(Y)+H_\alpha(Z) - H_\alpha(Y,Z)$ is strongly influenced by the scaling of $Y$ and $Z$.
We can eliminate the effect of the marginals by transforming the data so that all marginals become identical.

We achieve this by applying the rank-transform \cite{Con12,3Spe04}, together with centering, to the input data. The marginals of the rank-transformed data are discrete representations of the uniform distribution $U[0,1]$. That is, each transformed input variable attains the values $\frac{i-1/2}{N}$, for $N$ the number of datapoints and $i$ the index of the datapoints. Because of this transform, the R\'{e}nyi entropy of the transformed input variables becomes a straightforward quantifier of dependence. Since the transformation is monotonic in each of the dimensions, the ordering of the datapoints along each coordinate axis is preserved. The only change is that after the transformation, the distances between datapoints along the coordinate axes are equal (original distances are distorted). Furthermore, outliers are in general less pronounced in the rank-transform than in the original form. Altogether, the rank-transform leaves the structure of the original data intact. We give some numerical examples of the rank-transform and its effects towards the end of this section.

\subsection{Estimator of the R\'{e}nyi entropy}
It is not straightforward how to estimate the entropy of a given dataset. One approach is to estimate the distribution underlying the dataset, and compute the entropy from this estimated distribution. Both parametric and non-parametric methods (e.g. kernel density estimation, binning) can be used to estimate the distribution, however the results are quite sensitive to the details of the method used \cite{Silverman}. Therefore we employ a different approach for entropy estimation in this study, one
in which the data is used directly and no estimation of the distribution is needed so that we can circumvent the sensitivities due to estimating the distribution.
This approach can be used to estimate the R\'{e}nyi entropy with $0<\alpha<1$, however not the Shannon entropy. Hence our focus on the R\'{e}nyi entropy here.

In \cite{2Dud81}, Shannon entropy has been used for a uniformity test on the unit interval, and \cite{2Kra04} proposed two classes of estimators for the mutual information based on $k$-nearest neighbor distances. However, it has been shown \cite{2Gao15} that the number of samples required scales exponentially with the MI itself, which makes accurate estimation between strongly dependent variables almost impossible. An improved estimator is proposed, but requires the setting of an extra parameter.

Hero \& Michel describe in three papers \cite{2Her03,2Her02,2Her98} direct methods to estimate the R\'{e}nyi entropy, which are based on constructing minimal graphs spanning the datapoints in the domain. The estimator, which is asymptotically unbiased, is
\begin{equation}\label{eq:Renyiestimator}
\hat{H}_\alpha (X) = \frac{1}{1-\alpha} \left(\log \left(\frac{L_\gamma(X)}{N^\alpha}\right) - \log \beta_{L,\gamma}\right) = \frac{1}{1-\alpha} \log \left(\frac{L_\gamma(X)}{\beta_{L,\gamma}N^\alpha}\right),
\end{equation}
with $X$ denoting the dataset, $\alpha=(d-\gamma)/d$ with $d$ the dimensionality of the domain of $X$ and $N$ the number of data points. $\beta_{L,\gamma}$ is a constant only depending on $\gamma$ and the definition of $L_\gamma$. Finally,
$L_\gamma(X)$ is a functional, defined as
\begin{equation}\label{eq:gamma}
L_\gamma(X) = \min_{T(X)} \sum_{e\in T(X)} |e|^\gamma,
\end{equation}
where $T(X)$ is the set of spanning trees on $X$ and $e$ denotes an edge. The parameter $\gamma$ can be freely chosen within the interval $(0,d)$, where $d$ denotes the dimension (number of input variables) of $X$. Furthermore, $d$ must satisfy the constraint $d\geq 2$. The choice of $\gamma$ determines $\alpha$, and the given bound on $\gamma$ guarantees $\alpha$ to be in $(0,1)$. Estimator (\ref{eq:Renyiestimator}) is strongly consistent for $\alpha\in(0,1)$. Cases with $\alpha>1$ are not considered. If $\gamma$ is chosen to be $1$ and $|e|$ denotes the Euclidean distance between datapoints in $X$, then $L_\gamma(X)$ describes the length of the minimum spanning tree (MST) on the dataset. If the dependence between two input variables is computed, then $d=2$ and thus $\alpha=1/2$. The effect of varying $\gamma$ and thereby $\alpha$ is largely outside of the scope of this study. Note that $\gamma=1$ is computationally a convenient choice. Furthermore, it has been proven that for a related quantity, the $\alpha$-divergence, the theoretically optimal value for distinguishing between two densities which are close to each other on the basis of $\alpha$-divergence is $\alpha=1/2$ \cite{3Her01}. Also, the only value of $\alpha$ for which the $\alpha$-divergence is monotonically related to a true distance metric between the densities is $1/2$ \cite{3Her01}.

Here, we propose to use the R\'{e}nyi entropy with $\alpha=1/2$ instead of the Shannon entropy since (i) Equation (\ref{eq:Renyiestimator}) can not conveniently be computed for $\alpha=1$, because it requires taking the limit of $\gamma$ to 0, which also influences $\beta_{L,\gamma}$ and (ii) $\alpha=1/2$ is both numerically and theoretically a convenient choice as argued above.

\begin{figure}[ht!]
\centering
\subfloat[Independent data.]{\label{fig:1a}\includegraphics[width=0.45\textwidth]{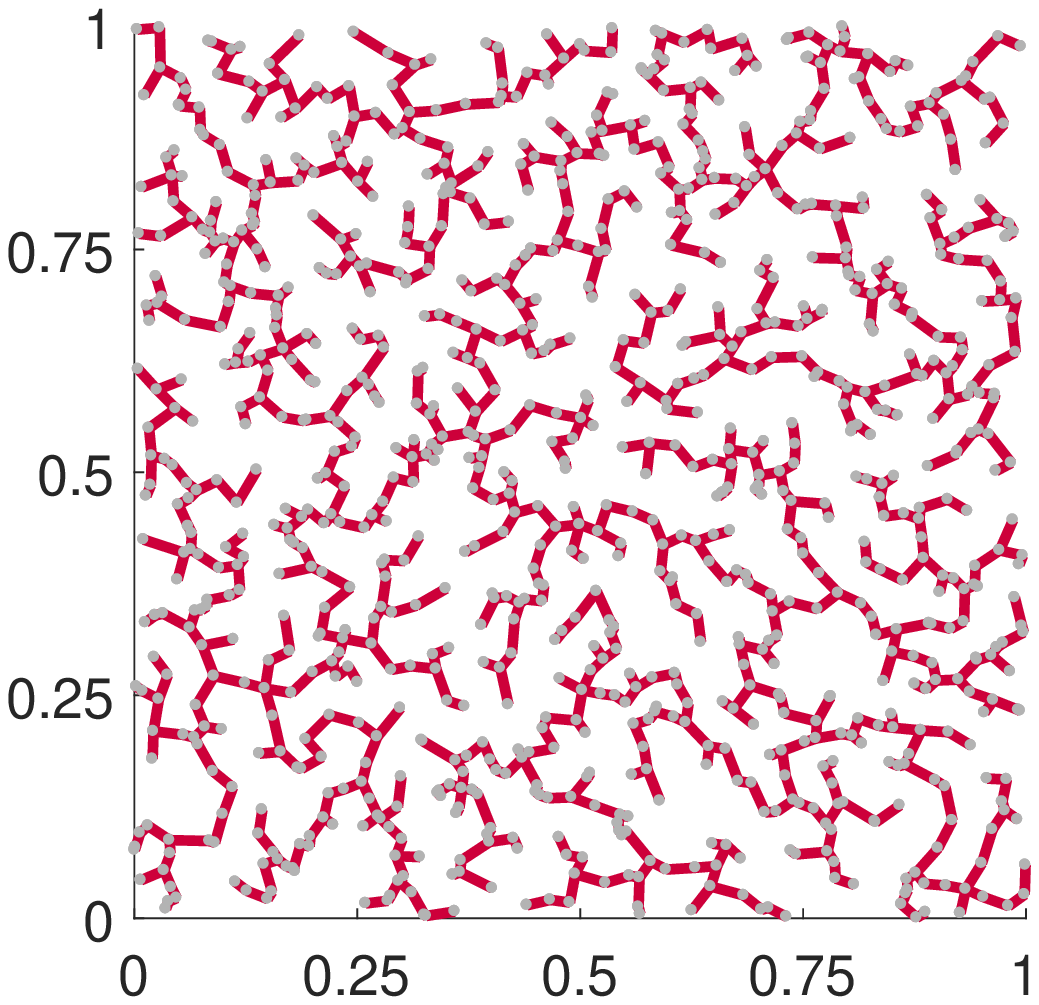}}\hspace{0.05\textwidth}
\subfloat[Dependent data.]{\label{fig:1b}\includegraphics[width=0.45\textwidth]{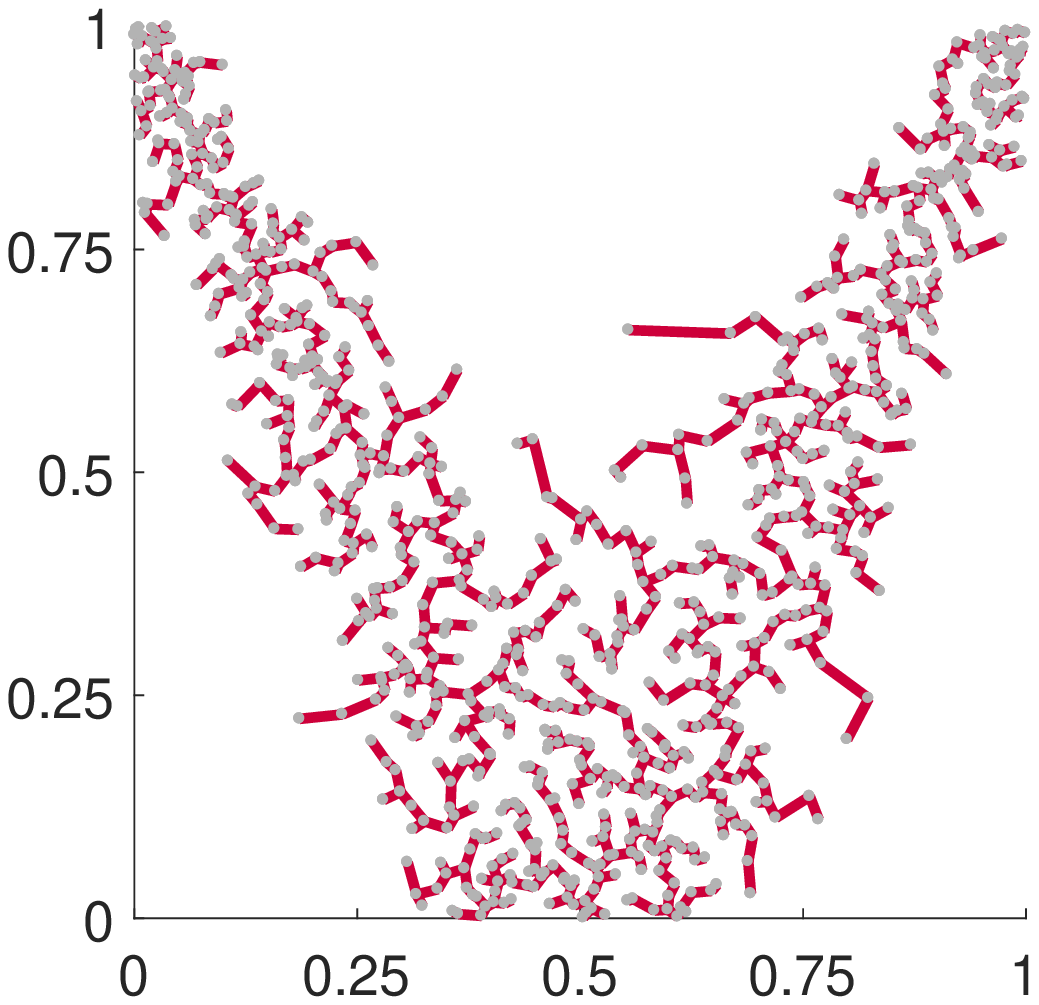}}
\caption{Illustration of the MST for two rank-transformed datasets. One dataset is sampled from a bivariate independent distribution, the other from a strongly nonlinear dependent distribution.}\label{fig:exMST} 
\end{figure}

Because of the scaling by the rank-transform as proposed earlier, in case of independence the MST will approximately cover the full domain (the unit square, if $d=2$). In case of dependence, either the density is nonuniform, or some areas of the domain contain no datapoints. In both cases, the average edge length decreases and so does the total length of the MST. This is illustrated in Figure \ref{fig:exMST}. If the length becomes shorter, the entropy estimate becomes smaller and the estimated dependence increases. Note that there is no assumption on the shape or structure of the dependence.

\subsection{Quantifier of dependence}
We propose to use the following quantity to estimate (or rather quantify)  the dependence between two (rank-transformed) variables:
\begin{equation}\label{eq:estimator}
H_\alpha^*(X) = \log\left(\frac{L_\gamma(X)}{N^\alpha}\right)
\end{equation}
where $X$ is a dataset of the variables containing $N$ datapoints. We set $\alpha=1/2$, as discussed in the previous section. With $\alpha$ and the dimension of the data ($d=2$) given, it follows that $\gamma=1$.
We note that $H_\alpha^*(X) = (1-\alpha)\hat{H}_\alpha(X) + \log(\beta_{L,\gamma})$. Thus, there is a linear relation between $H_\alpha^*(X)$ and $\hat{H}_\alpha(X)$, as $\alpha$ and $\beta_{L,\gamma}$ are constants as long as the dimension ($d$) and $\gamma$ do not change. Thus, for quantifying and comparing dependencies of different two-dimensional datasets, using (\ref{eq:estimator}) is equivalent to using (\ref{eq:Renyiestimator}), yet simpler in practice as we do not have to provide the constant $\beta_{L,\gamma}$, which cannot be obtained easily. In fact, the relation between (\ref{eq:estimator}) and (\ref{eq:Renyiestimator}) is affine.

Besides obtaining an ordering in terms of the dependency strength, we can use (\ref{eq:estimator}) to construct a reference level for distinguishing between independent and dependent variables.
This construction consists of computing $H_\alpha^*(X)$ on $r$ different datasets of rank-transformed bivariate uniformly distributed data. From this we obtain an empirical distribution, based on $r$ samples, for $H_\alpha^*(X)$ for independent variables (we recall that data from independent variables always leads to data uniformly distributed on a grid after the rank-transform). The reference level can then be defined as the $\eta$-quantile of the empirical distribution for $\eta$ small, since $H_\alpha^*(X)$ decreases with increasing dependence. This leads to the statistical test for dependence with the null hypothesis (independence) rejected if
\begin{equation}
H_\alpha^*(X) \leq \eta,
\end{equation}
in which $\eta$ is the 0.01 or 0.05 quantile of the empirical distribution for $H_\alpha^*(X)$ for independent variables.

\subsection{Proof of concept}\label{ssec:poc}
First, we compute the quantifier of dependence (\ref{eq:estimator}) on multiple datasets with varying distributions to analyze its behavior. Then, its convergence and robustness for increasing values of $N$ are studied. For the distributions considered in this section, it is possible to compute the R\'{e}nyi entropy with high accuracy for the scaled (i.e., rank-transformed) distribution such that a comparison can be made between the behavior of (\ref{eq:estimator}) and (\ref{eq:Renyi}).

\begin{figure}[ht!]
\centering
\subfloat[Original data (corner distribution).]{\label{fig:5a}\includegraphics[width=0.45\textwidth]{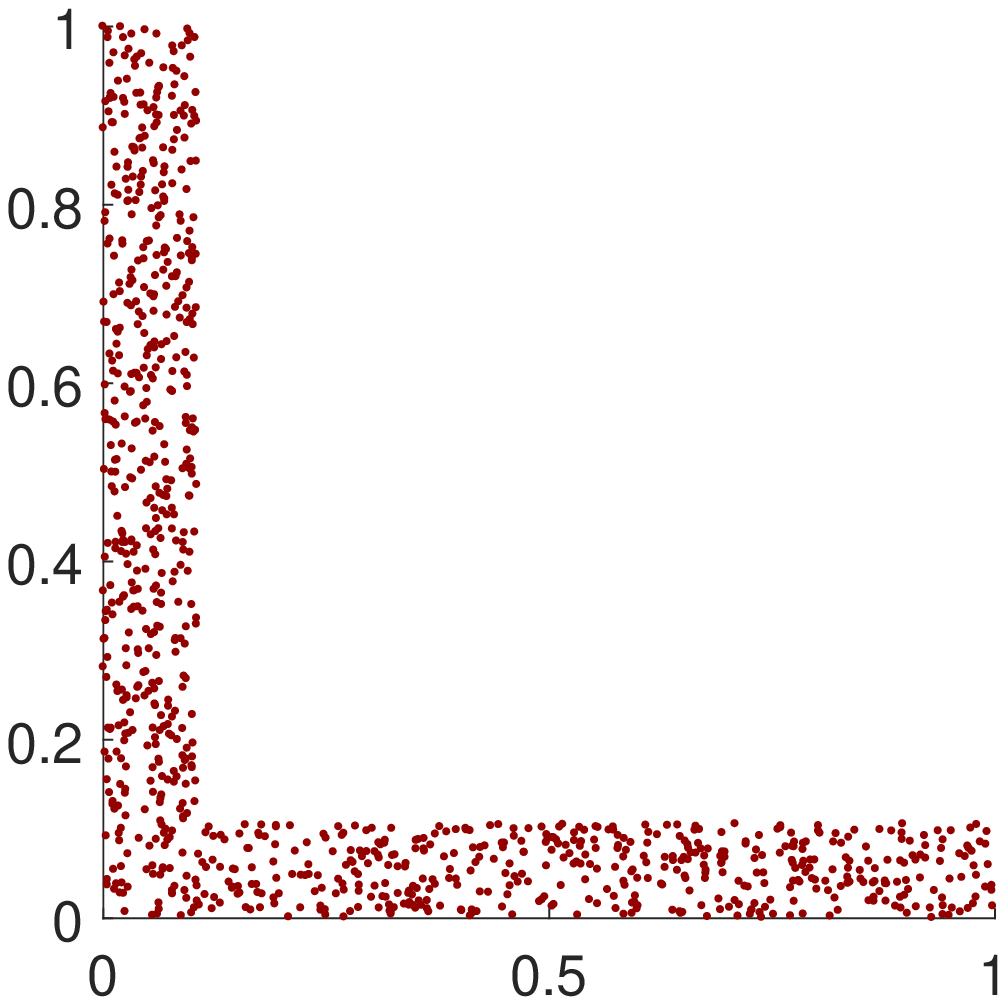}}
\subfloat[Rank-transformed data (corner distribution).]{\label{fig:5b}\includegraphics[width=0.45\textwidth]{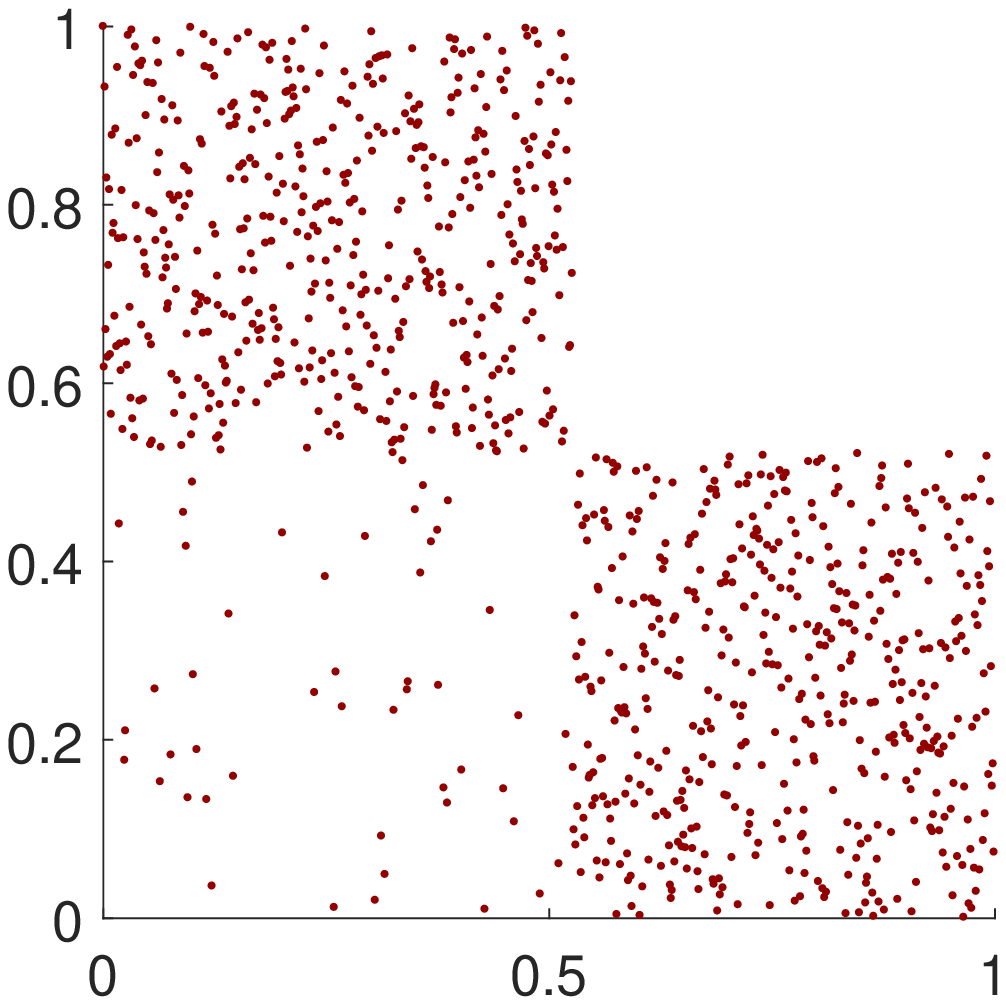}}\\
\subfloat[Original data (line distribution).]{\label{fig:5c}\includegraphics[width=0.45\textwidth]{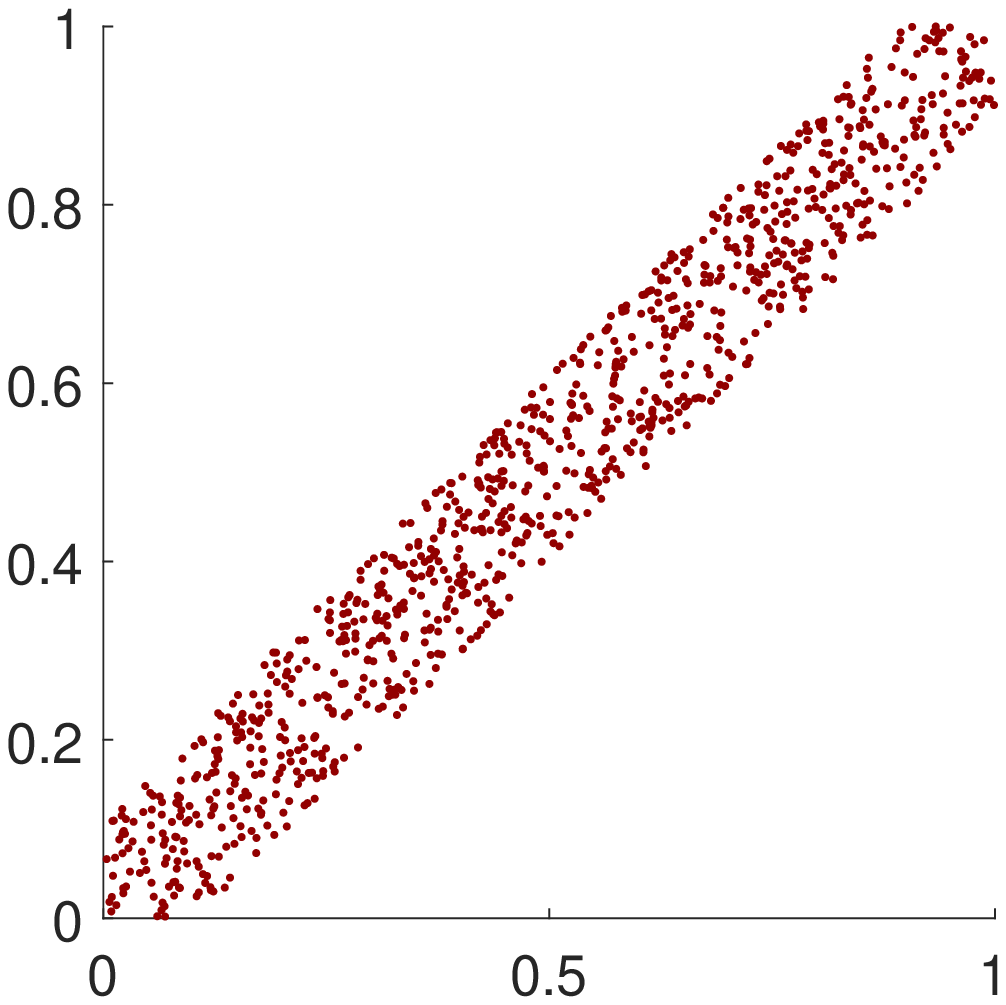}}
\subfloat[Rank-transformed data (line distribution).]{\label{fig:5d}\includegraphics[width=0.45\textwidth]{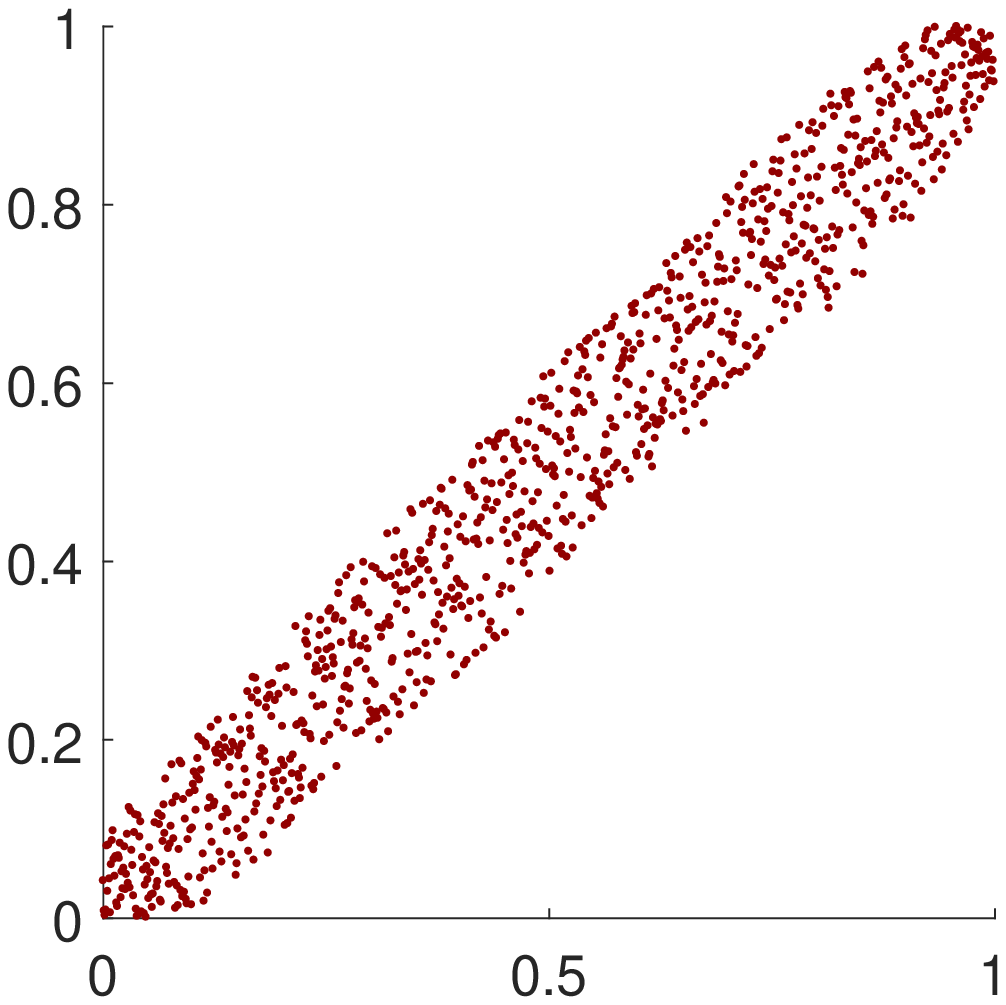}}
\caption{Example of the data and its rank-transform for the shape distribution with $A=0.8$.}\label{fig:Lscaling}
\end{figure}

We use datasets sampled from the following distributions: (i) a bivariate uniform distribution, (ii) a standard normal distribution with varying correlation coefficient $\rho$, (iii) a constant density on a region within the unit square with area $1-A$, in which the data is focused in the lower left corner (corner distribution), and (iv) a constant density on a region within the unit square with area $1-A$, but in which the data is focused on a symmetry axis (line distribution). For the latter two cases, the region includes values near the minimum and maximum of both variables, such that the ranges of the region for varying $A$ stay the same. These distributions are also referred to as shape distributions. Figure \ref{fig:Lscaling} shows examples of these distributions. For each of the four distributions, $r=10^2$ datasets are generated  with $N=10^3$ datapoints in every dataset. $\rho$ and $A$ are varied in steps of $0.10$ from $0.05$ to $0.95$.

\begin{figure}[ht!]
\centering
\includegraphics[width=0.6\textwidth]{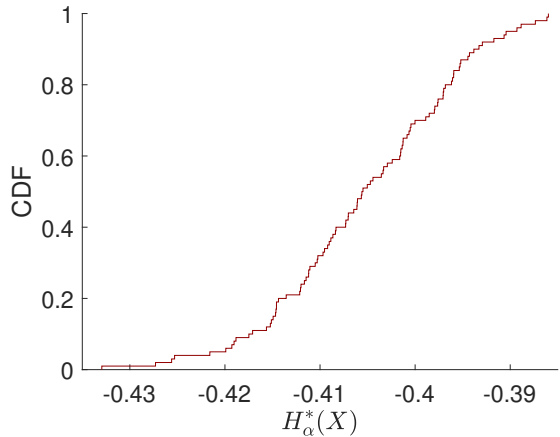}
\caption{Empirical distribution of $H_\alpha^*(X)$ (\ref{eq:estimator}) for the bivariate uniform distribution.}\label{fig:unif}
\end{figure}

In Figure \ref{fig:unif}, the empirical cumulative distribution function (CDF) of $H_\alpha^*(X)$ (\ref{eq:estimator}) is plotted for data from the bivariate uniform distribution. It can be seen that the distribution is concentrated in a rather narrow interval. At the end of this section we explore how this interval becomes more or less narrow as the amount of data $(N)$ changes.

\begin{figure}[ht!]
\centering
\subfloat[Normal distribution.]{\label{fig:3a}\includegraphics[width=0.33\textwidth]{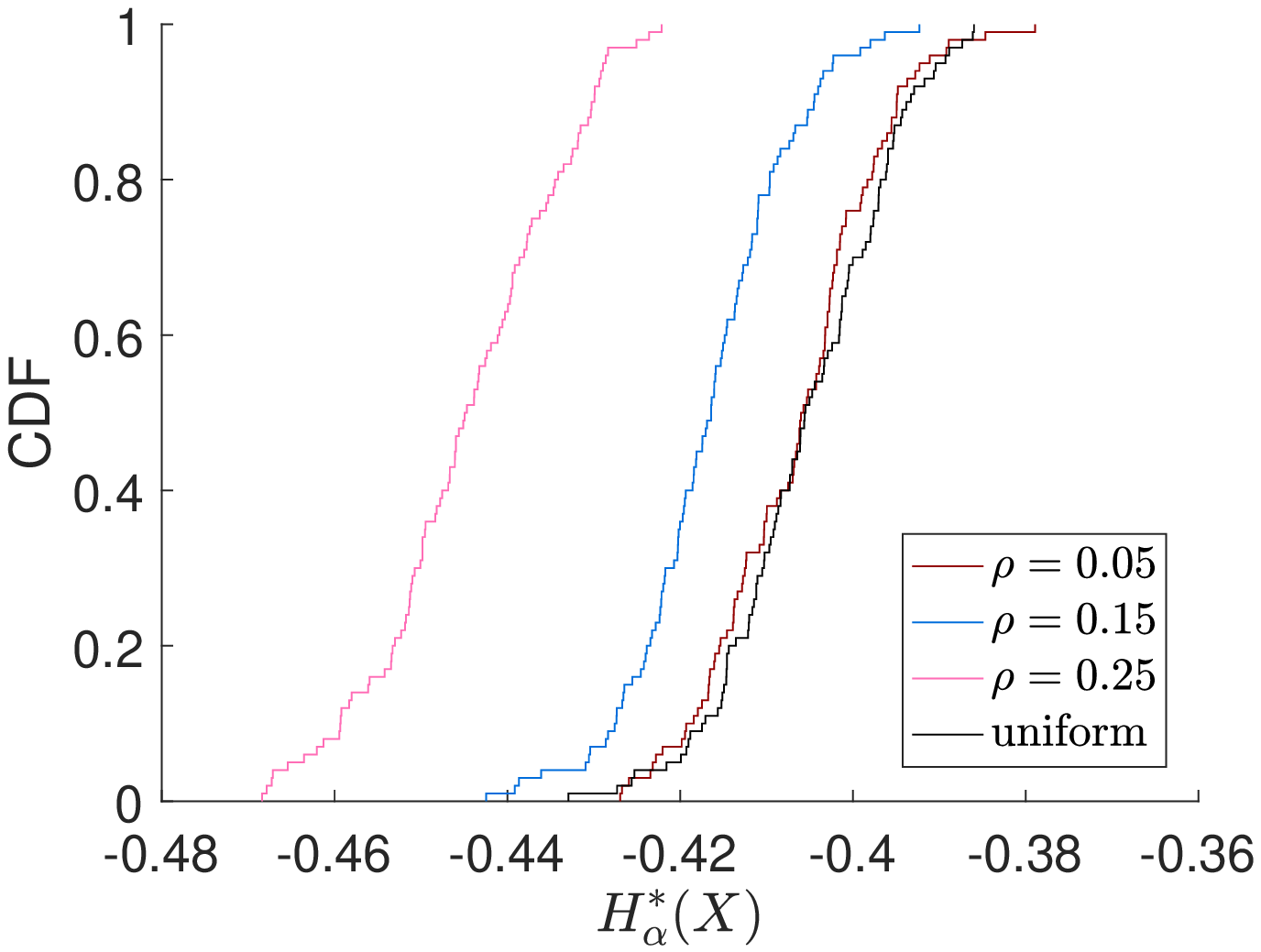}}
\subfloat[Corner distribution.]{\label{fig:3b}\includegraphics[width=0.33\textwidth]{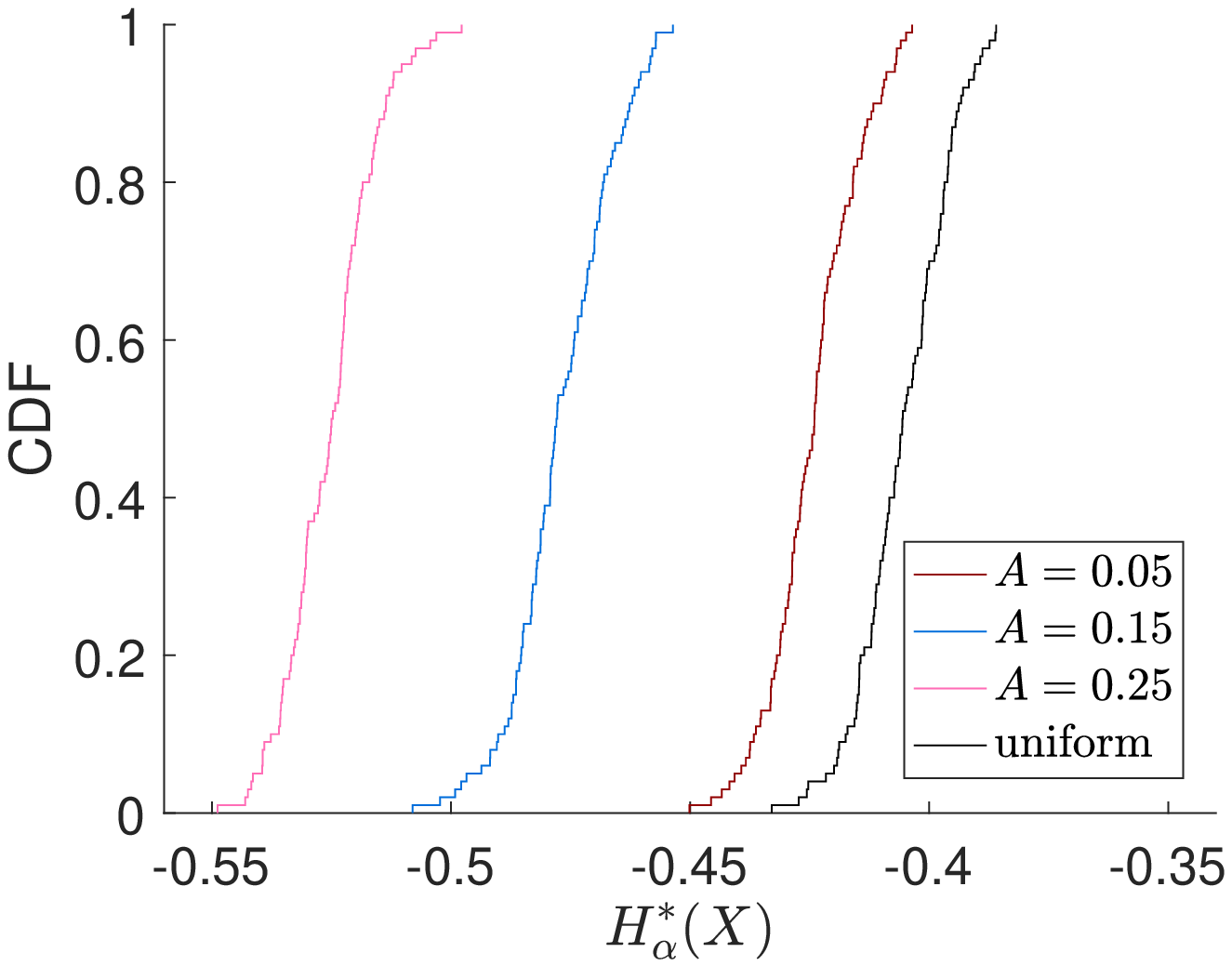}}
\subfloat[Line distribution.]{\label{fig:3c}\includegraphics[width=0.33\textwidth]{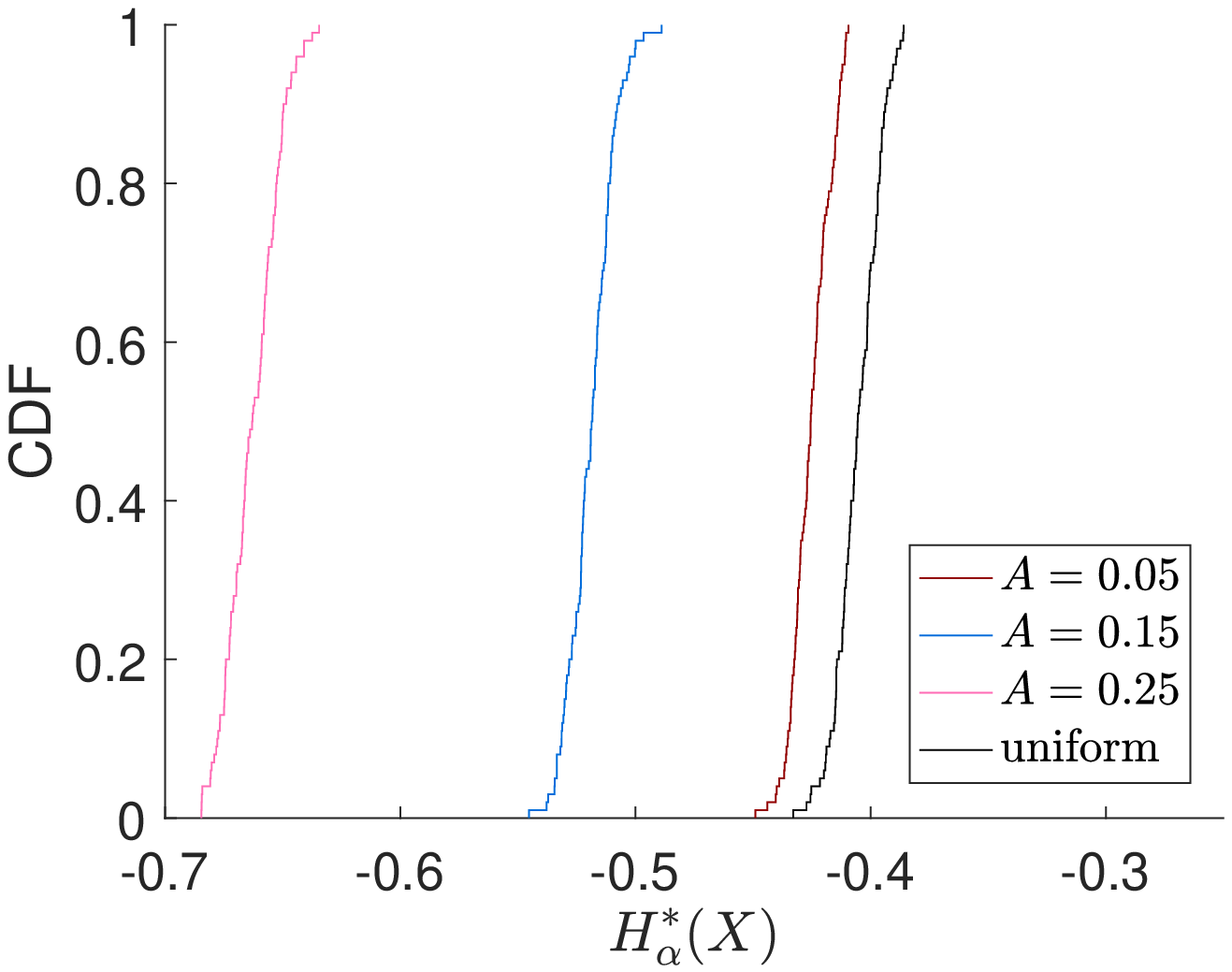}}
\caption{Empirical distributions of $H_\alpha^*(X)$ for data from the normal distribution and the two shape distributions, for small $\rho$ and $A$.}\label{fig:other}
\end{figure}

In Figure \ref{fig:other}, the empirical CDF of $H_\alpha^*(X)$ is plotted for both the normal distribution (case (ii)) and the two shape distributions (case (iii) and (iv)), for different small values of $\rho$ and $A$. More precisely, for each dataset sampled from one of these distributions, we first apply the rank-transform and then evaluate (\ref{eq:estimator}) for the rank-transformed data.
The empirical CDF for the uniform distribution, already shown in Figure \ref{fig:unif}, is included in black for comparison. If $\rho=0$ or $A=0$, there is no dependence anymore and the empirical CDF obtained from the normal distribution or shape distributions coincides with the empirical CDF from the bivariate uniform distribution (modulo differences due to finite sample size, $N < \infty$).

\begin{figure}[ht!]
\centering
\subfloat[R\'{e}nyi entropy.]{\label{fig:4a}\includegraphics[width=0.45\textwidth]{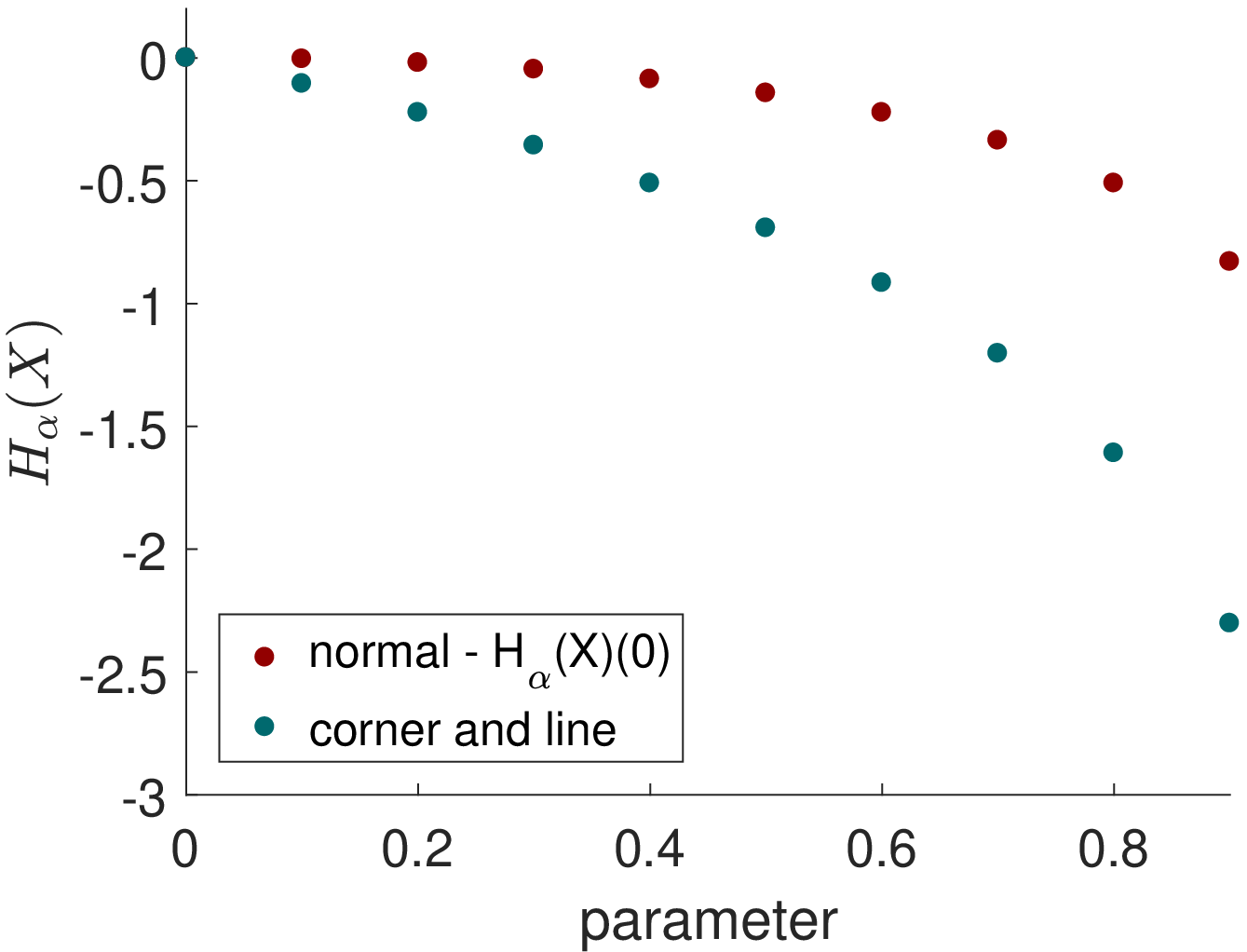}}\hspace{0.05\textwidth}
\subfloat[Quantifier of dependence.]{\label{fig:4b}\includegraphics[width=0.45\textwidth]{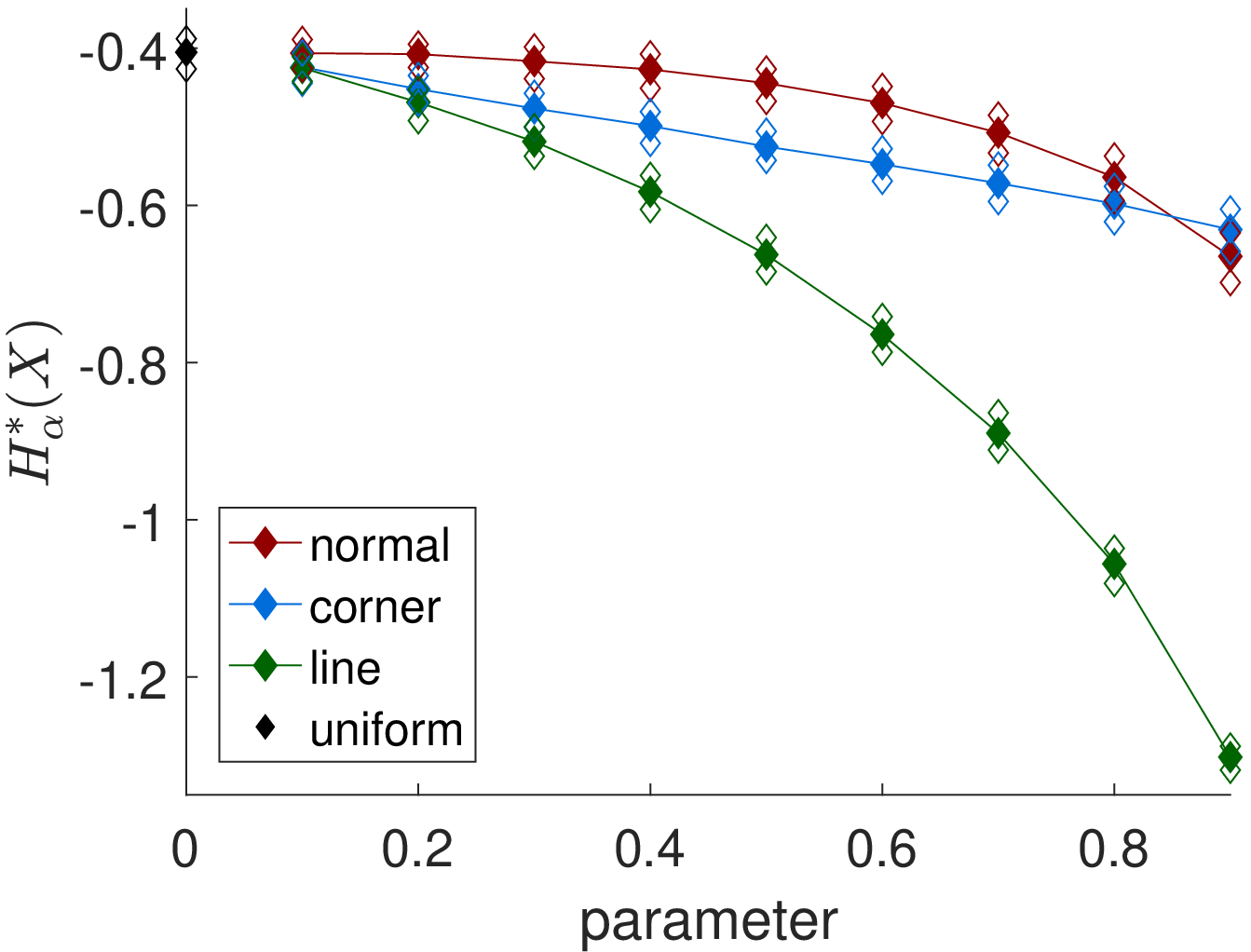}}
\caption{Comparison of the R\'{e}nyi entropy (left) and quantifier of dependence (right) for the normal distribution and shape distributions, with varying parameters ($\rho$ and $A$). The R\'{e}nyi entropy is computed exactly using (\ref{eq:Renyi})
 without transforming the data. For visualization purposes, the values for the normal distribution have been translated by its value for $\rho=0$, which is $2\log(2\sqrt{2\pi})$.
The quantifier (\ref{eq:estimator}) is evaluated using data that is sampled from the distribution and then rank-transformed. We show the mean and 95\% confidence intervals of $H_\alpha^*(X)$. Note that the numerical values of the entropy and $H_\alpha^*(X)$ are not supposed to coincide, cf. (\ref{eq:Renyiestimator}) and (\ref{eq:estimator}).} \label{fig:CIs}
\end{figure}

It can be seen that for the shape distributions, the quantifier of dependence is more distinctive for small values of $A$ than it is for the normal distribution with small $\rho$. Hence, weak dependencies in the shape distributions can more easily be detected than in the normal distribution. To demonstrate the behavior of $H_\alpha^*(X)$ for the full range of values of $\rho$ and $A$, we plot the mean together with the empirical 95\% confidence intervals in Figure \ref{fig:4b}. In the case of $\rho=A=0$, the distribution of the uniform (independent) distribution is recovered, indicated by a different color.

From Figure \ref{fig:4a}, it can be seen why the differences between the CDFs are small in case of the normal distribution: for the normal distribution, the entropy is nearly flat as a function of $\rho$, for small $\rho$ values. Figure \ref{fig:4b} shows the rank-transform has a more severe effect on the corner than on the line distribution, since the estimates for the corner distribution decrease slower. The distortion of the shape of the graph for the corner distribution between Figure \ref{fig:4a} and \ref{fig:4b} is caused by the distortion of the distances between data points after the rank-transform. Note that both shape distributions have the same R\'{e}nyi entropy if the rank-transform is not applied. The effect of the rank-transform in this case can be seen in Figure \ref{fig:Lscaling}. The closer $A$ is to 1, the more the data falls in the two boxes for the corner distribution, due to the combination of skewness and discontinuity. This does not hold for the line distribution. In the limit of $A\rightarrow 1$ and $N\rightarrow\infty$, the entropy of the rank-transformed corner distribution goes to $-\log(2)$, while it goes to $-\infty$ for the rank-transformed line distribution. Hence, it is consistent that $H_\alpha^*(X)$ does not go to $-\infty$ for the corner distribution. We note that the shape distributions are quite artificial and have discontinuous density, while datasets in practice are usually samples from a distribution with a smooth density and a less artificial shape.

\begin{figure}[ht!]
\centering
\includegraphics[width=0.6\textwidth]{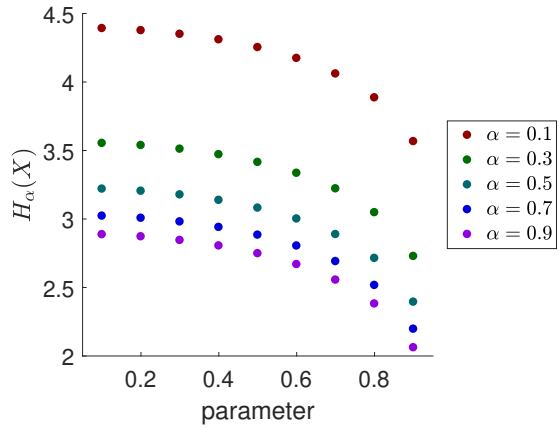}
\caption{R\'{e}nyi entropy for the normal distribution with parameter $\rho$ (correlation coefficient), for varying $\alpha$.}\label{fig:varyalpha}
\end{figure}

For completeness, we give the values for $H_\alpha(X)$ with $X$ the normal distribution with parameter $\rho$ and for varying $\alpha$ in Figure \ref{fig:varyalpha}. It can be seen that the differences are small, apart from a shift that is constant in $\rho$. This constant shift has no consequence for our goal of quantifying dependencies by means of entropy. In particular, for all values of $\alpha$ shown, the entropy is nearly flat as a function of $\rho$ as $\rho \to 0$. Thus, other values than $\alpha=1/2$ offer no advantage in this respect, supporting our choice for $\alpha=1/2$.

\begin{figure}[ht!]
\centering
\includegraphics[width=0.6\textwidth]{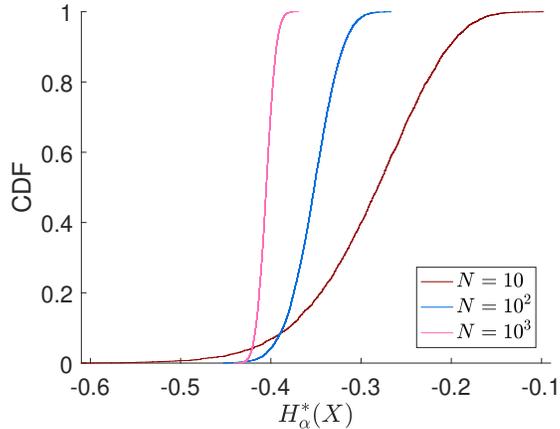}
\caption{Empirical CDF of $H_\alpha^*(X)$ for the bivariate uniform distribution. The CDF becomes narrower for larger datasets (increasing $N$).}\label{fig:unifNvary}
\end{figure}

We conclude this section by investigating the effect of varying the size of the dataset, $N$. We compute the empirical CDF of $H_\alpha^*(X)$ using $N=10$, $10^2$ and $10^3$. For this computation we increase $r$ to $10^4$ to obtain a smooth CDF. The resulting empirical CDFs are shown in Figure \ref{fig:unifNvary}. It can be seen that the distribution becomes narrower with increasing $N$, as was to be expected. Thus, larger $N$ makes the comparison of estimates easier due to the smaller confidence intervals involved.

We note that for higher values of $N$, the computations become expensive due to the computation of the edge lengths before computing the MST. The order of this operation is $O(N^2)$. Kruskal's method for computing MSTs \cite{2Kru56} runs in $O(N^2\log(N))$ time, but the steps are faster than the ones needed for computing the edge lengths. Prim's method for computing MSTs \cite{2Pri57} can be faster ($O(N^2+N\log(N))$), but the implementation is more involved. The high cost of computing the exact MST motivates us to investigate approximate algorithms in the next section.

\section{Approximation methods}\label{sec:methods}
For large datatsets (i.e., large $N$), the computational cost of evaluating the estimator (\ref{eq:estimator}) becomes prohibitively high, due to the computational complexity of constructing the MST. In this section we discuss methods to reduce the computational cost by approximating the value of (\ref{eq:estimator}) evaluated on a large dataset. We consider three types of approximation: first, MSTs can be computed on multiple subsets of the data. The second type aggregates data points into clusters, such that only one MST calculation is performed on the cluster centers. The third type clusters the data points, constructs MSTs on each cluster and combines them in a smart way \cite{2Zho15}.

\subsection{Sampling-based MST}
In this method, the dataset is split in $K$ subsets and $H_\alpha^*(X)$ (\ref{eq:estimator}) is computed on each of the subsets. Thus, $K$ estimates of $H_\alpha^*(X)$ are obtained. Their arithmetic mean becomes the new estimate, while their variance is a measure for the quality of the new estimate. The number $K$ can be chosen freely, although it represents a trade-off between accuracy and computation time (since both accuracy and computational time decrease with increasing $K$).

The splitting can be done in various ways, of which random and stratified are the most straightforward ones. In the random splitting, datapoints are allocated randomly to one of the $K$ subsets, which all have size $N/K$ (rounding is neglected). In the proportional splitting, $k$ clusters are generated with the $K$-means method \cite{3Cel13,Ste56} and these are considered the strata. The data points in each stratum are then proportionally allocated to the $K$ subsets. The number of clusters involved ($k$) would be another parameter, which we choose equal to $N/K$.

\subsection{Cluster-based MST}
Another way of reducing the computational burden would be to cluster the data in a large number of clusters and compute the MST on the cluster centers. In this case, the clustering method can be chosen, as well as the number of clusters. To be consistent with the previous method, we define here the number of clusters to be $k=N/K$, such that the number of points in the MST is similar to the previous method. Furthermore, it is possible to include the size (weight) of the clusters as well. Two different clustering methods are investigated: $K$-means \cite{Ste56} and PCA-based clustering \cite{2Egg18}. The construction of the MST is performed both without and with weighting. The weighting is harmonic, which implies that in the construction of the MST, the edge between datapoints $i$ and $j$, $|e_{ij}|$, is replaced by
\begin{equation}
W_{ij}|e_{ij}|, \quad W_{ij} = \frac{2k}{\frac{1}{w_i}+\frac{1}{w_j}},
\end{equation}
where $k$ is the number of cluster and $w_i$, $w_j$ are the weights of the clusters, computed by the fraction of datapoints in that cluster.

\subsection{Multilevel MST}
This method is developed by Zhong et al.\ \cite{2Zho15} and is also called FMST (fast MST). It is based on the idea that to find a neighbor of a data point, it is not necessary to consider the complete dataset of size $N$. In the FMST method, the dataset is partitioned in $\sqrt{N}$ clusters via the $K$-means clustering method, and MSTs are constructed on each of the clusters. In a next step, a MST is constructed on the cluster centers to determine which clusters get connected at their boundaries. Between these clusters, the shortest connecting edges are computed. Because this is heuristic, the complete process is repeated with the midpoints of the edges connecting different MSTs being the initial cluster centers for the $K$-means method. The resulting two MSTs are merged into a new graph and its MST is computed as final outcome of the method. \cite{2Zho15} report good results with this method, which has a computational complexity of $O(N^{1.5})$. Errors occur only if points that should be connected end up in different clusters twice, which does not occur often and has only a minor effect. For high-dimensional datasets, erroneous edges are included more often, but the effect thereof is smaller. Both are due to the curse of dimension.

\subsection{Comparison}
\begin{figure}[ht!]
\centering
\subfloat[Independent data.]{\label{fig:7a}\includegraphics[width=0.45\textwidth]{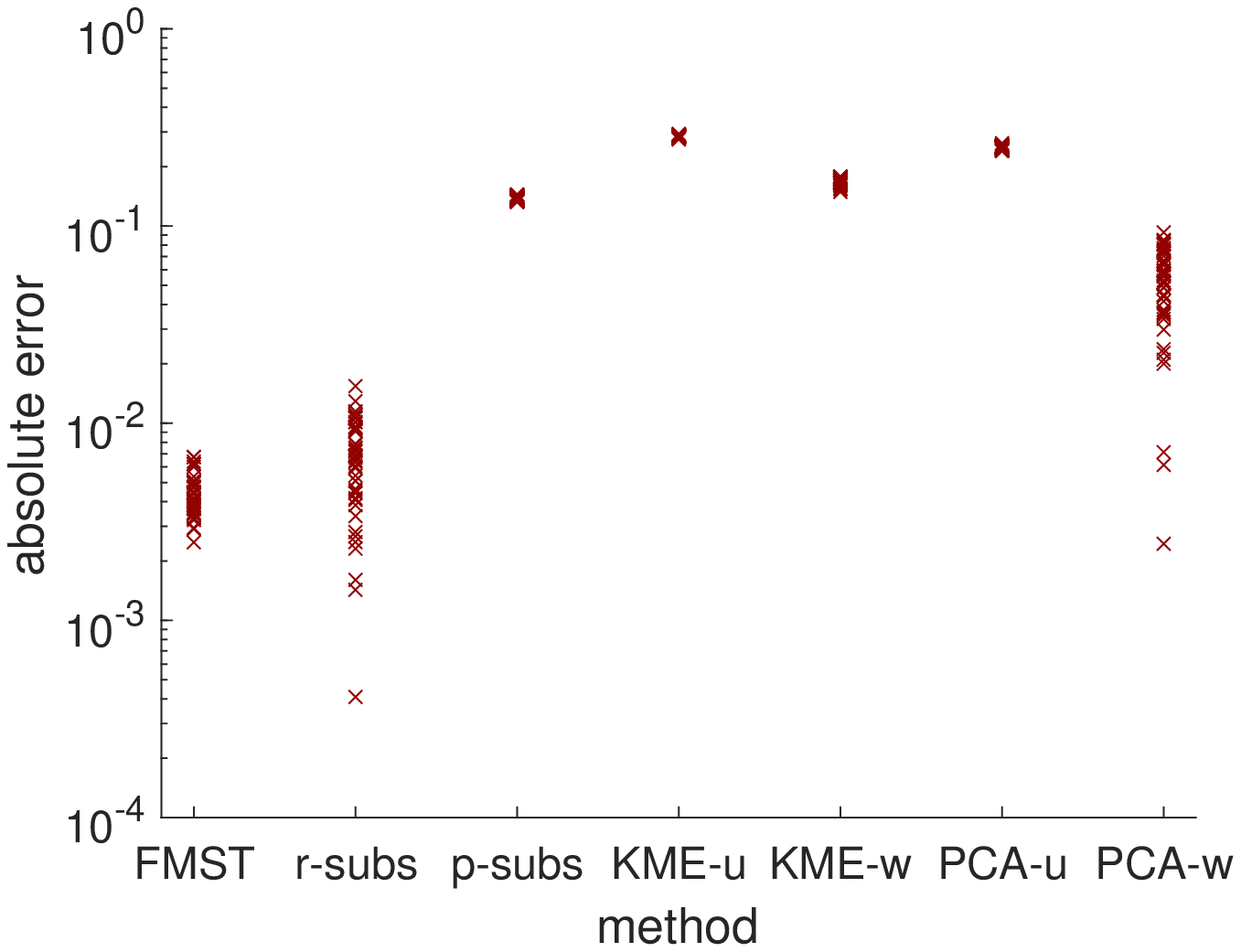}}
\subfloat[Dependent data.]{\label{fig:7b}\includegraphics[width=0.45\textwidth]{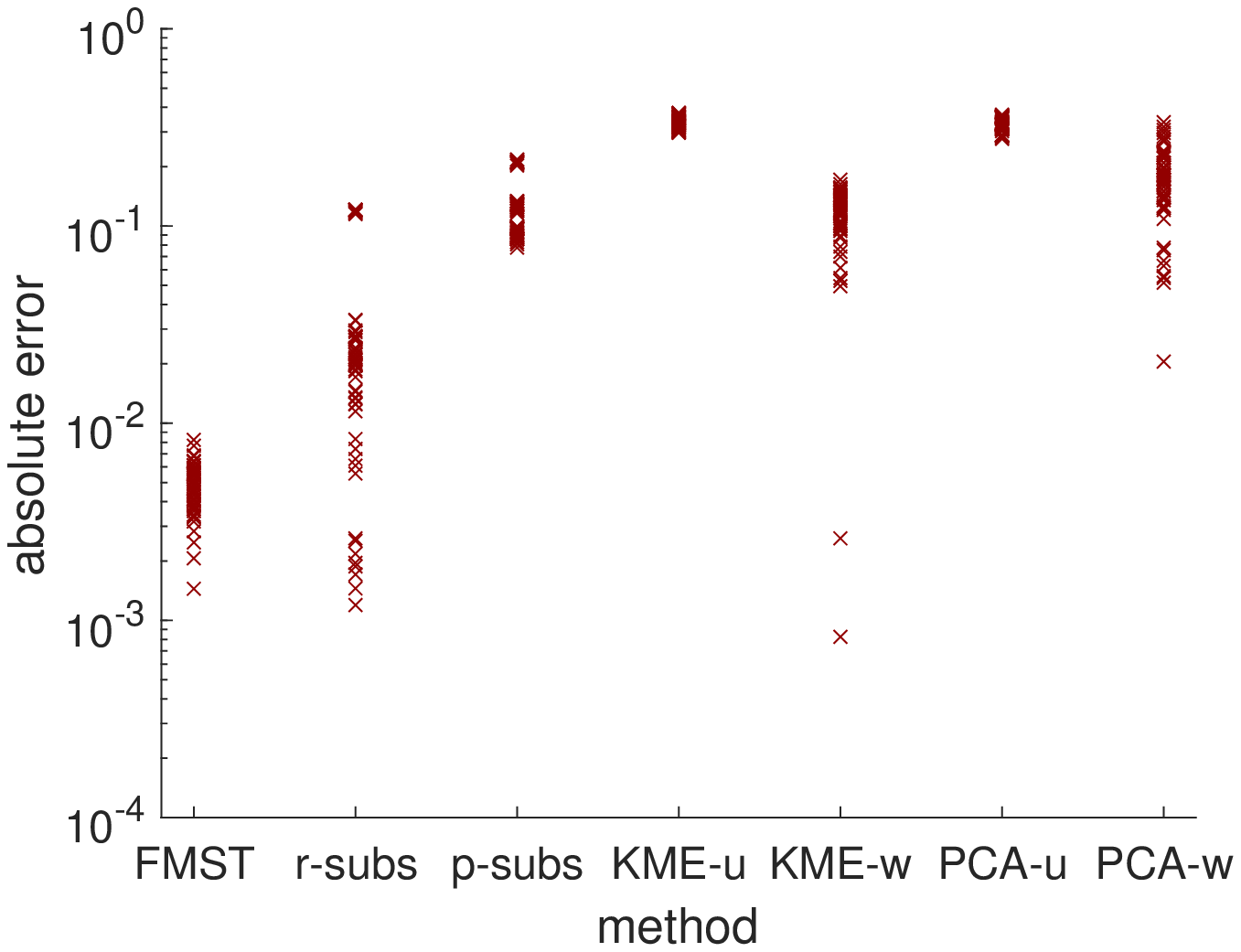}}
\caption{Error estimates for different approximation methods. The figure on the left is for independent data, while the figure on the right is made with dependent data. From left to right, the methods are: FMST, random sampling from subsets-based MST, stratified sampling from subsets-based MST, $K$-means cluster-based MST (unweighted and weighted) and PCA-based cluster-based MST (unweighted and weighted). }\label{fig:sec3NLL}
\end{figure}

The accuracy, robustness and computational time of the methods explained before are tested on $r=50$ bivariate uniform datasets with $N=10^4$ and $r=60$ bivariate strongly dependent datasets with $N=10^4$. The strongly dependent datasets are chosen as projections of the dataset used in Section \ref{ssec:Ishi}. Hence, the dataset exists of combinations of $x$, $y$, $z$ and $I(x,y,z)$. Since this leads to 6 projections per dataset, we choose $r=60$ in this case. In these cases, it is computationally expensive but still feasible to compute the full (exact) MST, and we do so to be able to compare results from the three approximate methods with the method based on the full MST. For the approximate methods, we choose the parameter $K$ to be $10$. The results can be found in Figure \ref{fig:sec3NLL}. We mention here the effect of the implementation of $K$-means, as the choice of initialization can affect the results. We use the $K$-means++ initialization algorithm, repeat the clustering 10 times from different initializations and select the result with the best clustering criterion.
The error estimate is computed as the absolute difference between the approximated and exact value of $H_\alpha^*(X)$. The FMST has consistently small error, while other methods have larger error.
In our opinion, it is important the approximation is accurate and robust, i.e., has little variation in its error over multiple runs of data with the same distribution. Furthermore, it should perform well for both dependent and independent data. Therefore, we propose to use the FMST to compute approximations to the MST in case where the dataset is large (say, $N\geq 10^4$).

\section{Validation of the proposed FMST estimator}\label{sec:results}
In this section, we further investigate using the FMST method with the estimator (\ref{eq:estimator}).
First, the effect of the size of the dataset $N$ is investigated together with the robustness of the FMST estimator for relatively small $N$. Then, the FMST estimator for dependence is tested for behavior and consistency using datasets sampled from the three distributions considered before (uniform, normal and shape distributions).

\begin{figure}[ht!]
\centering
\includegraphics[width=0.6\textwidth]{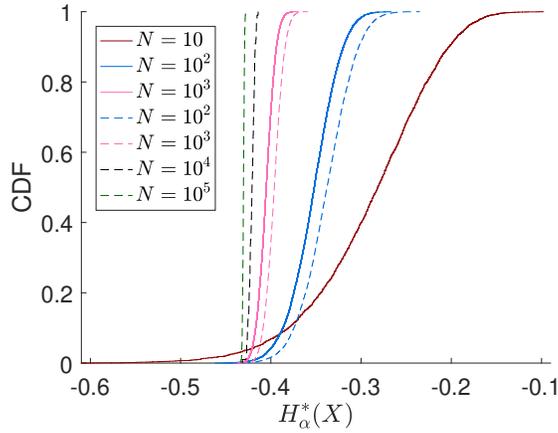}
\caption{Empirical distribution for the uniform distribution for varying $N$. The solid lines refer to the distributions based on the MST, while the dashed lines refer to the distributions based on the FMST. Results using MST are limited to $N\leq 10^3$ because of high computational cost. }\label{fig:unifNvaryFMST}
\end{figure}

It is straightforward that approximations to the MST using the correct edge distances overestimate the length of the MST. Therefore, we compare the empirical distributions of $H_\alpha^*(X)$ based on the MST (Figure \ref{fig:unifNvary}) to the ones based on the FMST in Figure \ref{fig:unifNvaryFMST}. The number of repetitions is $r=10^4$ for $N=10^2$, $10^3$ and $r=10^2$ for $N=10^4$ and $N=10^5$. It can be seen that the distributions are different, but only shifted. Furthermore, the width of the distribution decreases to almost zero for $N=10^4$ and $N=10^5$. This means the estimator is robust to sampling effects. Note that there is still a difference between the distributions for $N=10^4$ and $N=10^5$. This is due to the bias caused by approximating the MST, which reduces for increasing $N$.

\begin{figure}[ht!]
\centering
\subfloat[Mean and confidence intervals for FMST.]{\label{fig:9a}\includegraphics[width=0.45\textwidth]{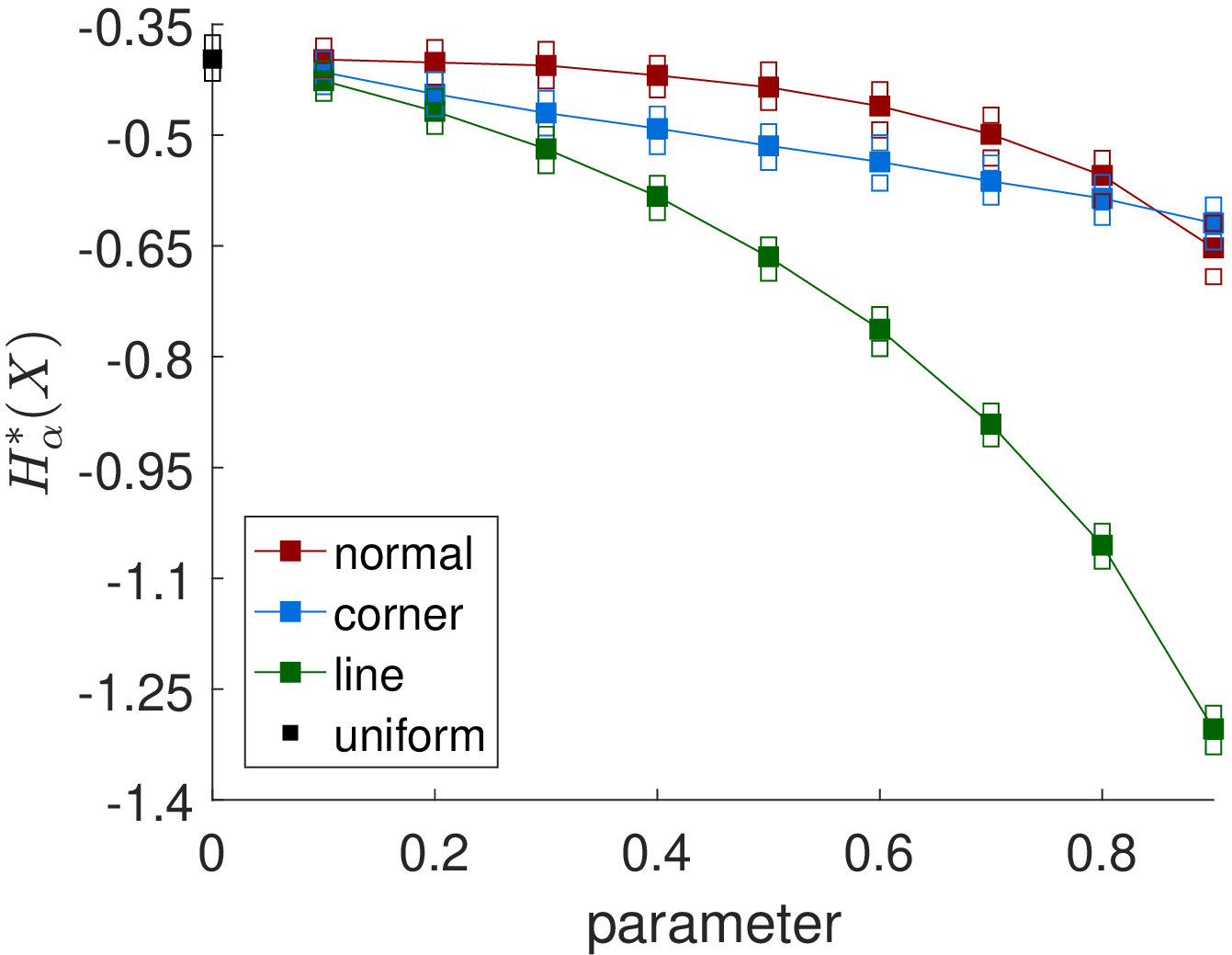}}
\subfloat[Comparison with MST.]{\label{fig:9b}\includegraphics[width=0.45\textwidth]{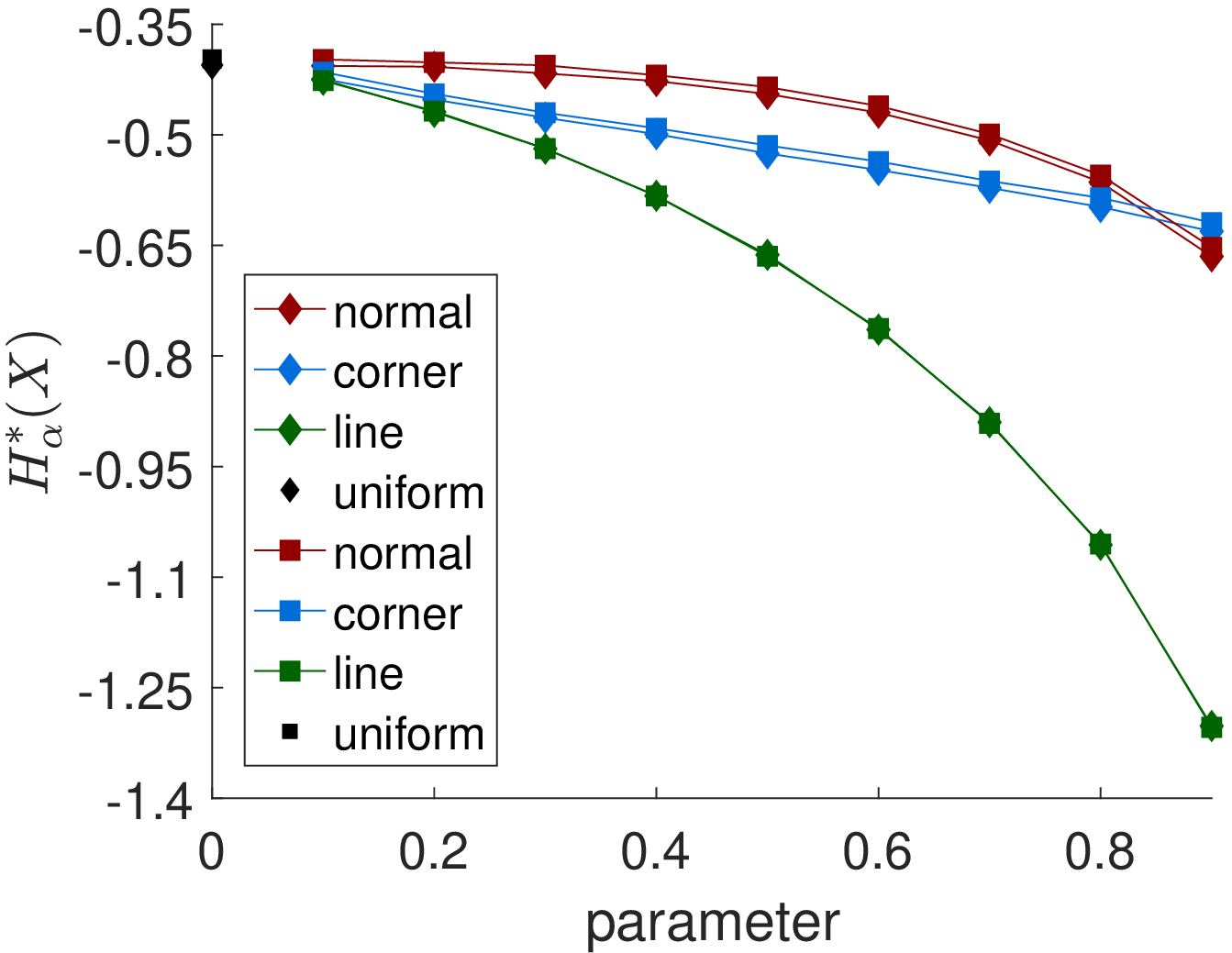}}
\caption{Behavior of the FMST estimator for two different types of data distribution. On the left the mean and empirical 95\%-confidence intervals, while on the right the means of the MST (diamonds) and FMST (squares) estimates can be compared (see also Figure \ref{fig:CIs}).}\label{fig:CIsFMST}
\end{figure}

We repeated the experiment from Section \ref{ssec:poc}, but now with the estimator based on FMST to evaluate the effect of the FMST approximation on different data distributions. The results are in Figure \ref{fig:CIsFMST}. Again, the behavior of the estimator based on FMST is the same as for MST, but only shifted upwards. This bias is very small for the line distribution. In these tests, $N=10^3$ in order to compare them to the previous results.

Based on these results, we conclude that the FMST estimator is a good and robust approximation of the MST estimator. It has a positive bias, but for the purpose of comparing and ranking strength of dependencies as quantified by the FMST estimator, this bias does not pose a problem.

\section{From dependency quantification to sensitivity analysis}\label{sec:5}
As already briefly mentioned in the introduction, the method proposed here to quantify dependencies via the MST can contribute to performing SA. If we have a $m+1$-dimensional dataset consisting of data for $m$ input variables $X_i$ ($i=1, ..., m$) and one output variable $Y$, we can quantify the dependencies by computing $H_\alpha^*(X_i,Y)$ for all $m$ combinations $(X_i,Y)$. The result can help to quantify sensitivity. In particular, it enables us to rank the input variables from most important (strongest dependency of input and output, largest sensitivity) to least important (weakest dependency of input and output, smallest sensitivity).

We note that the first step in our proposed method consists of performing a rank-transform on the data, before computing the MST. Clearly, this step is needed as well in case the data includes an output variable.

The ranking of input variables from most to least important and quantification of sensitivity of the output with respect to different inputs is a main goal of SA. Different from other SA methods, there is no assumption of independent inputs needed for the method we propose. The idea to use the quantification of dependencies between inputs and outputs by means of entropy for purposes of SA was proposed before by \cite{2Aud08,2Liu06}, however different methods were used in these studies to estimate entropy (or a proxy thereof) compared to what we propose. Furthermore, in these studies the case with multiple dependent inputs was not considered.

It must be mentioned here that by considering combinations of one input and one output, higher-order effects (interactions) cannot be explored. MSTs on higher-dimensional combinations of input and output variables may be a way to generalize this approach for higher-order effects. We leave this for further study. Note that the number of combinations to be computed grows fast with the number of inputs.

For comparison purposes, in the next section we include an example with independent input variables and one output variable, of which the Sobol's (total) indices can be computed analytically. These serve as a reference to test the proposed method on its qualities for sensitivity analysis. We compare our estimates to Sobol's and Sobol's total indices.

\section{Testcases}\label{sec:tests}
The proposed methods are applied to two testcases in this section. In the first case, we use the Ishigami function  \cite{3Ish90} and evaluate it using randomly sampled synthetic data as inputs. We consider data from two different input distributions, one without dependencies (uniform) and one with strong dependencies. On all combinations of variables (input and output), the dependence is quantified and, for the uniform dataset, compared to the values of the Sobol (total) indices.
In the second testcase, measurements on wave conditions at sea are investigated, of which the resulting sediment transport is computed.

\subsection{Ishigami function}\label{ssec:Ishi}
This test function is from Ishigami \& Homma \cite{3Ish90} and its Sobol (total) indices \cite{3Sob01} can be computed analytically. The test function is given by
\begin{equation}
I(x,y,z|a,b) = (a+bZ^4)\sin(X) + a\sin^2(Y),
\end{equation}
where
\begin{equation}
X = -\pi+2\pi x, \quad Y = -\pi+2\pi y, \quad Z = -\pi + 2\pi z.
\end{equation}
The parameters $a$ and $b$ are chosen to be $7$ and $0.1$, respectively, in accordance with \cite{2Cre09}. One dataset is four-dimensional and uniformly distributed, while the other one is generated as
\begin{equation}\label{eq:polydata}
x \sim U(-2,2), \quad y \sim x^2 + \frac{1}{2}N(0,1), \quad z \sim x^3 + \frac{1}{2}N(0,1), \quad u \sim x^4 + \frac{1}{2}N(0,1),
\end{equation}
and range-normalized to the unit hypercube. Note that there is one variable ($u$) which is not used in the Ishigami function. This is on purpose to show the effect of confounders. The MSTs are approximated with both the MST and the FMST method and (\ref{eq:estimator}) is used to compute the dependencies both between input variables and between the input variables and the output variable. The ordering is $(x,y)$, $(x,z)$, $(x,u)$, $(x,I)$ et cetera. Because of this ordering, the pairs numbered 4, 7, 9 and 10 represent combinations of one input and the output variable ($(x,I), (y,I), (z,I), (u,I)$, respectively).

\begin{figure}[ht!]
\centering
\subfloat[Uniform dataset.]{\label{fig:10a}\includegraphics[width=0.45\textwidth]{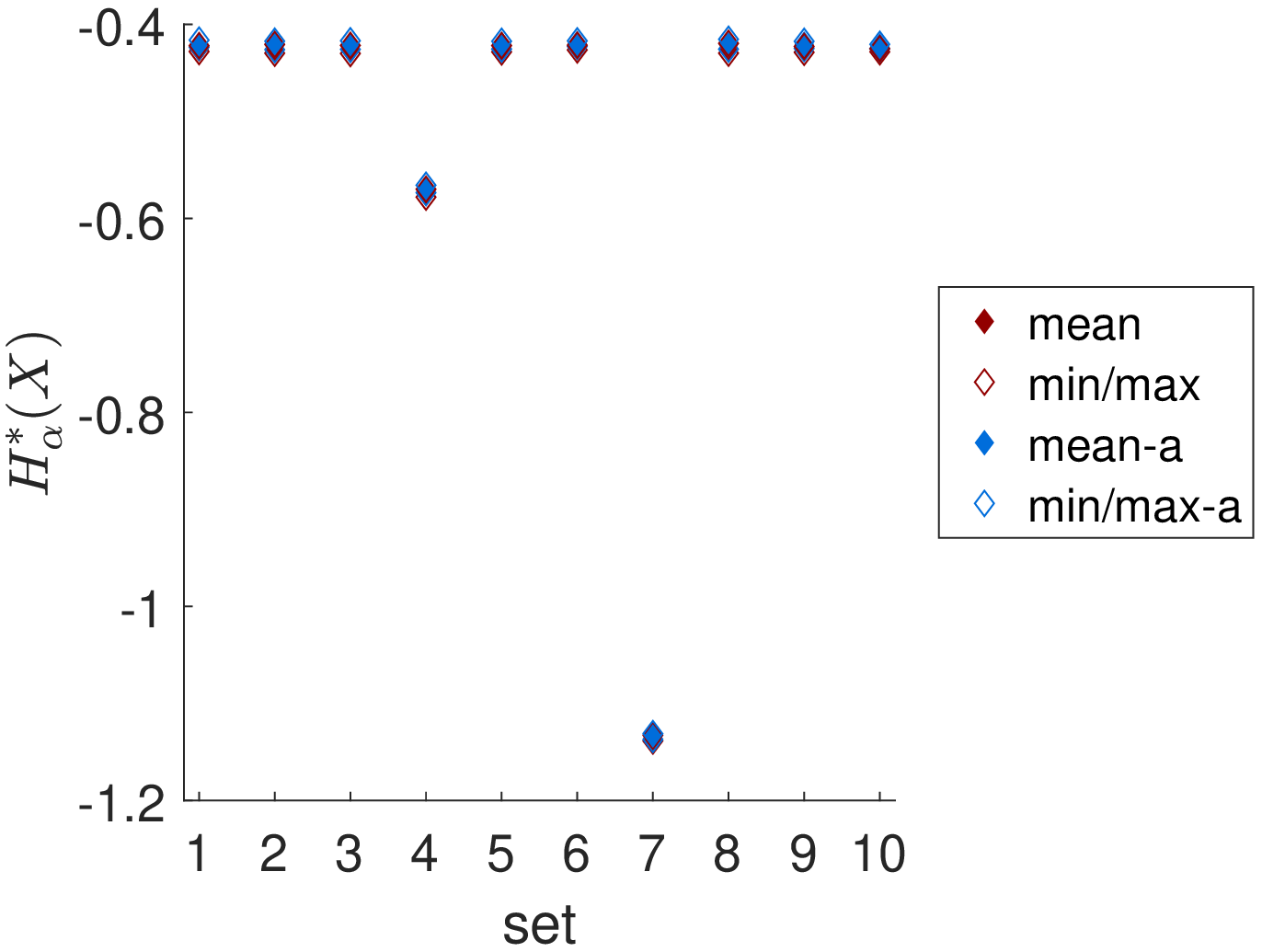}}\hspace{0.05\textwidth}
\subfloat[Dependent dataset.]{\label{fig:10b}\includegraphics[width=0.45\textwidth]{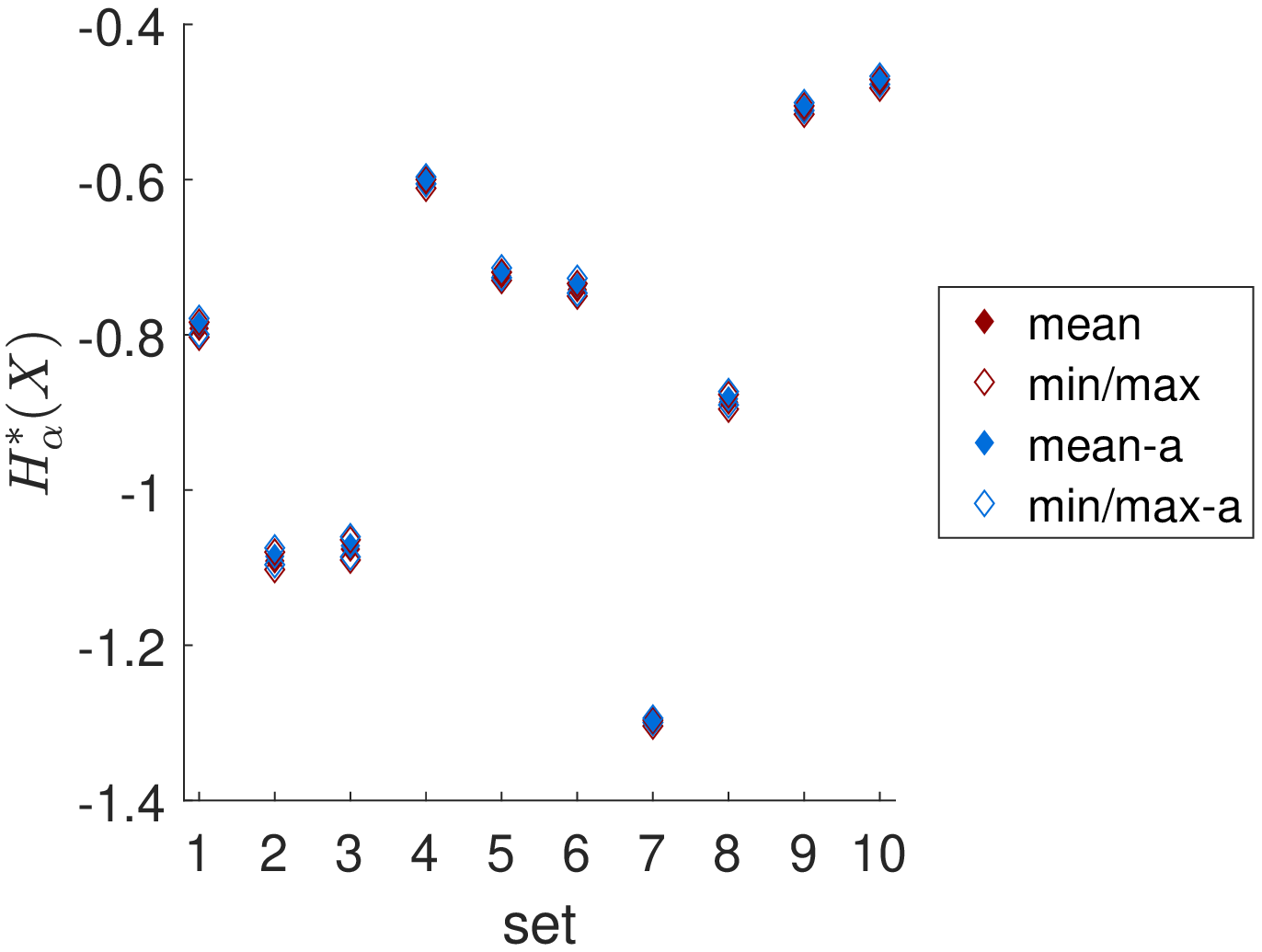}}
\caption{Estimates of $H_\alpha^*(X)$ for two datasets. Note the difference in the range of the $y$-axis. }\label{fig:resultsIshi}
\end{figure}
The datasets have been generated $r=10$ times with $N=10^4$ samples each and the results are in Figure \ref{fig:resultsIshi}.
First of all, it can be seen that the estimates are robust, as for each set (or pair) the estimates are all in a narrow interval (indicated by the minimum and maximum estimates).

Furthermore, it can be seen that for the uniform dataset (Figure \ref{fig:10a}), only the input variables $x$ and $y$  are found to have an effect on the output. For most sets, the dependence estimate is slightly below -0.4, the value attained in case of independence. One can compare this with Figure \ref{fig:unifNvaryFMST}, where it can be seen that in case of two uniformly distributed independent variables, the estimate is slightly below -0.4 in case of $N=10^4$. Indeed, the input variables are all independent in case of this uniform dataset. Also, by construction the Ishigami function does not depend on $u$. The independence of $z$ and $I$ is nontrivial, however it is consistent with the analysis using Sobol indices (as discussed below). $I$ is not directly dependent on $z$, although there is an interaction effect with $x$. This interaction effect is not detected with the analysis of pairwise dependencies.

\begin{figure}[ht!]
\centering
\includegraphics[width=0.35\textwidth,trim=15mm 0mm 2cm 5mm,clip]{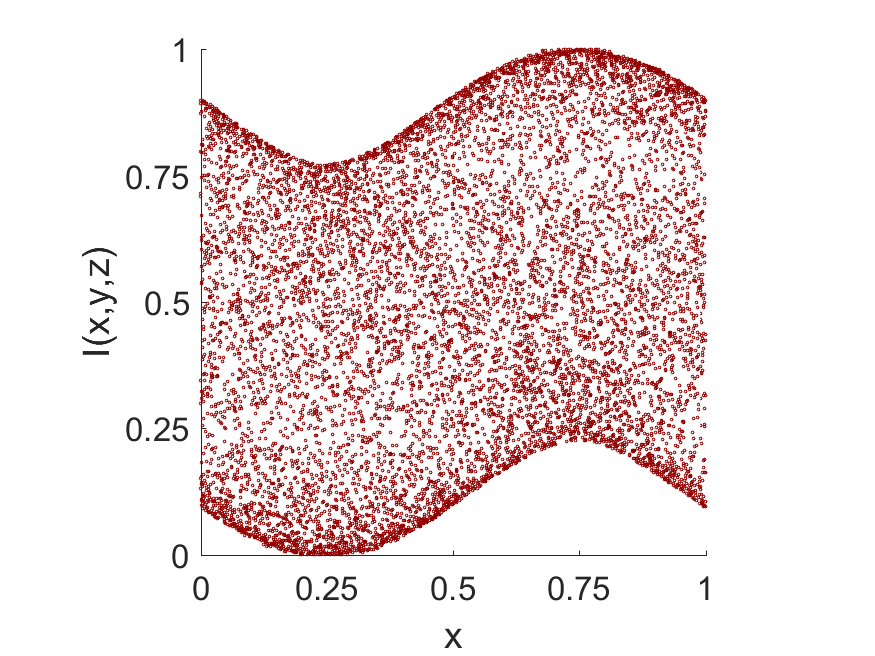}
\includegraphics[width=0.35\textwidth,trim=15mm 0mm 2cm 5mm,clip]{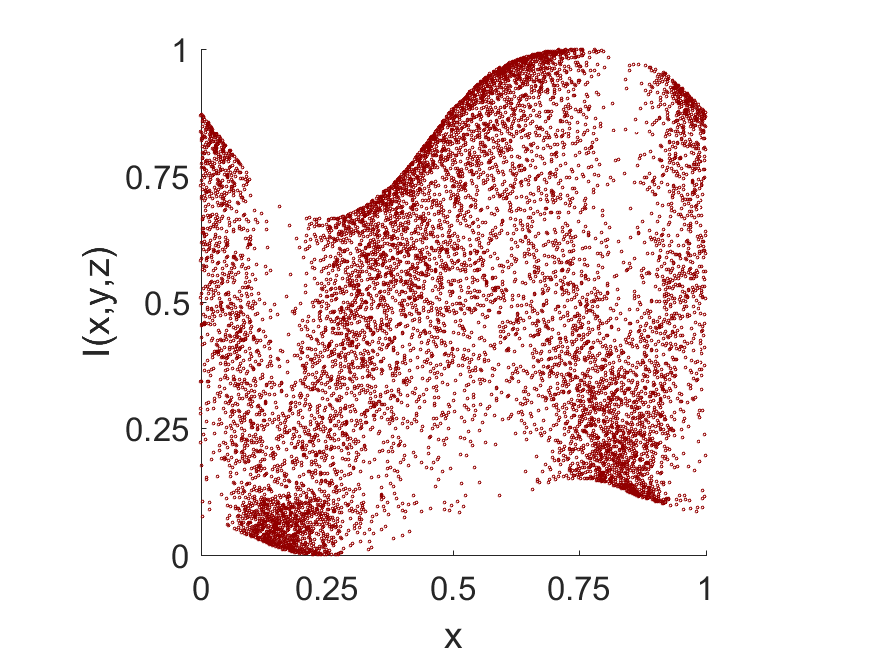}
\includegraphics[width=0.35\textwidth,trim=15mm 0mm 2cm 5mm,clip]{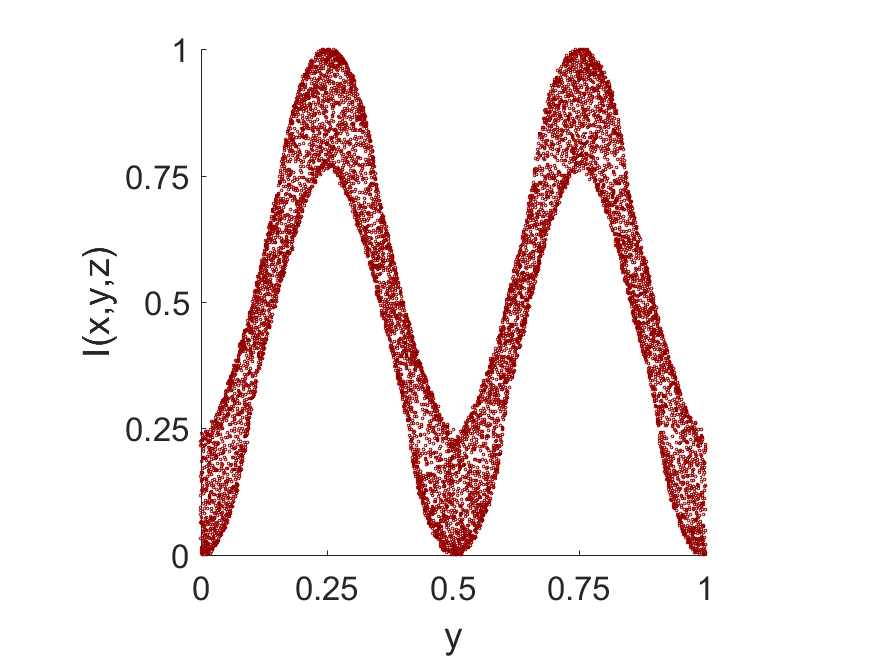}
\includegraphics[width=0.35\textwidth,trim=15mm 0mm 2cm 5mm,clip]{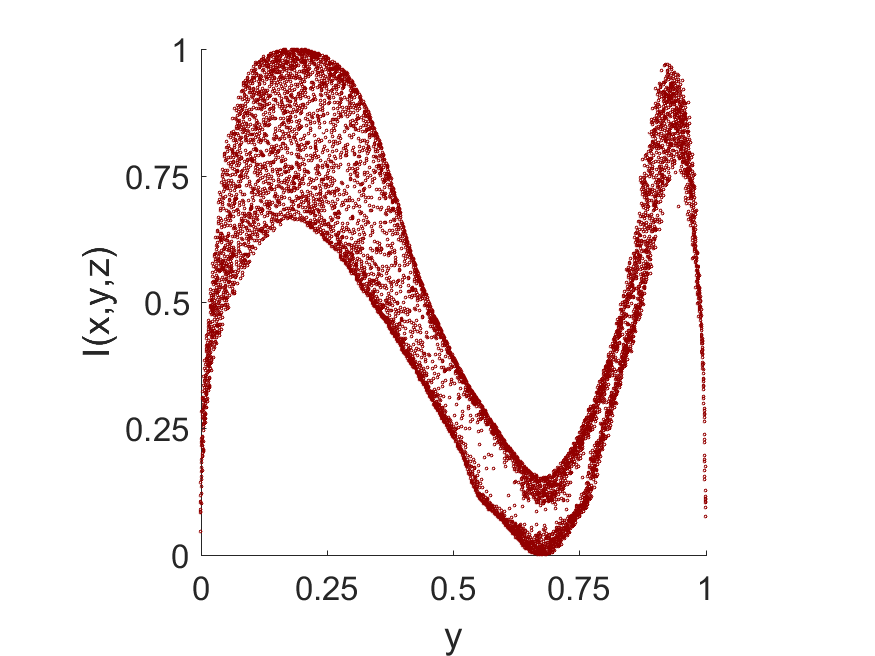}
\includegraphics[width=0.35\textwidth,trim=15mm 0mm 2cm 5mm,clip]{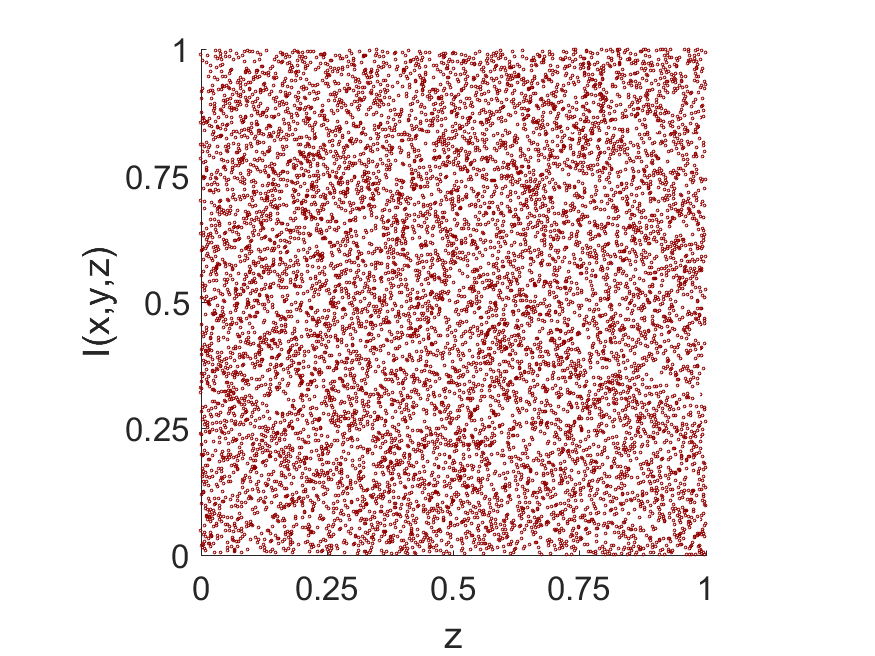}
\includegraphics[width=0.35\textwidth,trim=15mm 0mm 2cm 5mm,clip]{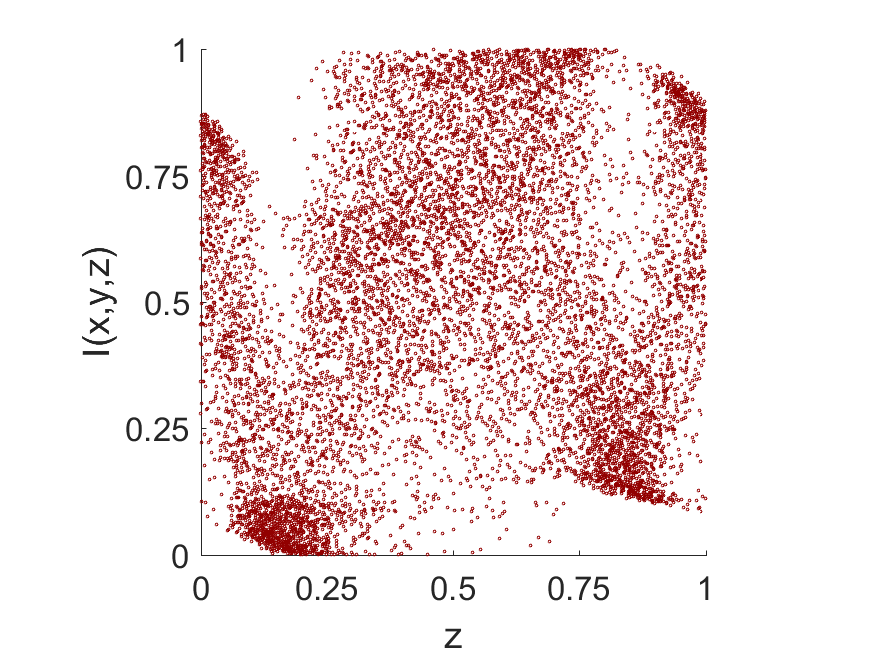}
\includegraphics[width=0.35\textwidth,trim=15mm 0mm 2cm 5mm,clip]{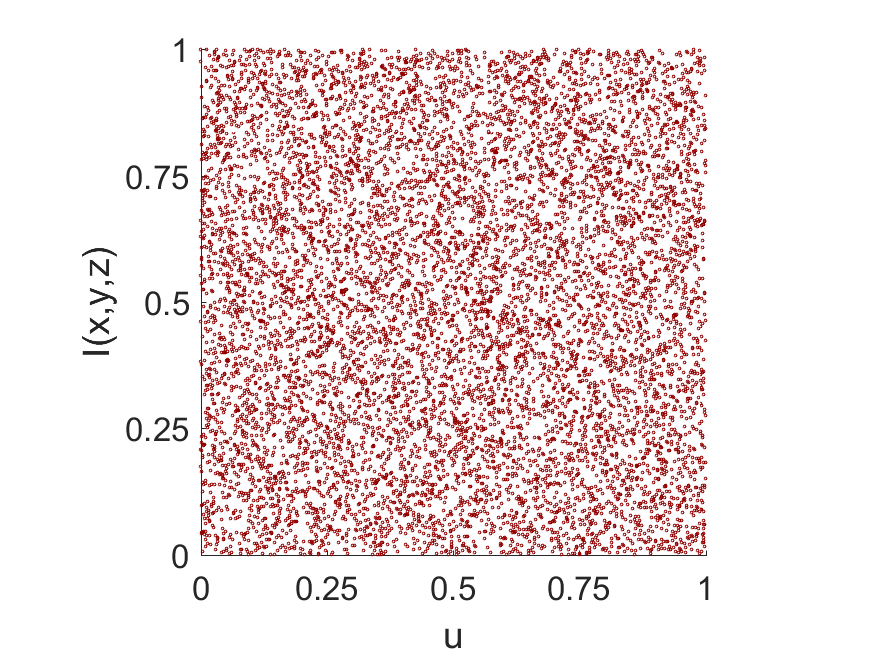}
\includegraphics[width=0.35\textwidth,trim=15mm 0mm 2cm 5mm,clip]{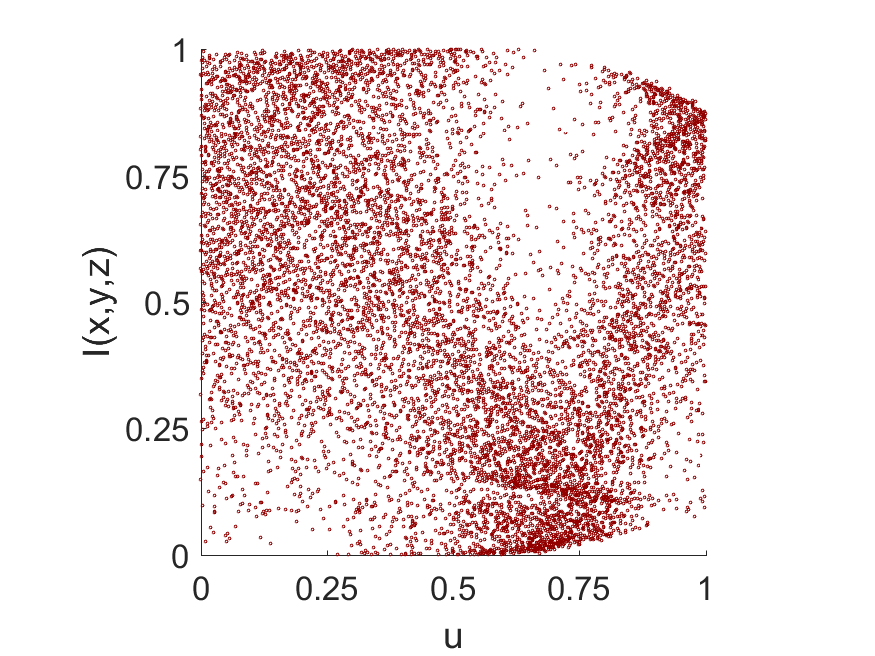}
\caption{Scatterplots for the Ishigami function. The dependent input data has a large effect on the output distribution. }\label{fig:Ishi}
\end{figure}

The exceptions to independence are sets 4 and 7 (pairs $(x,I)$ and $(y,I)$), for which the $H_\alpha^*$ estimate is lower, indicating dependence. The dependence of $y$ and $I$ is stronger than between $x$ and $I$, visible in the smaller value for set 7. The (in)dependence of $I$ on the various inputs is illustrated in Figure \ref{fig:Ishi} where scatterplots of the input-output pairs are shown.

For the strongly dependent dataset, the analysis is less straightforward. Again, the combinations $(x,I)$ and $(y,I)$ are strongly dependent, while $(z,I)$ shows a weak dependence. Because of the structure of the dataset, all input variables are mutually dependent. Furthermore, the dependencies within the input data are stronger than the dependencies of input variables with the output, except for the $(y,I)$ combination. The relative ordering of dependencies (strongest for $(y,I)$, weaker for $(x,I)$, very weak for $(z,I)$ and $(u,I)$) is consistent with the definition of the Ishigami function and the intuition gleaned from the right panels of Figure \ref{fig:Ishi}.

We continue by comparing the means of the estimates for the input-output combinations for the uniform dataset based on FMST with the values of its Sobol indices and Sobol total indices in Table \ref{tab:sobol}.
\begin{table}[ht!]
\caption{Comparison of the proposed estimator and Sobol's indices for the Ishigami function with independent (uniformly distributed) input variables.}
\label{tab:sobol}
\centering
\begin{tabular}{|r|ccc|}
\hline
& $S_i$ & $S_{Ti}$ & $H_\alpha^*$\\
\hline
$x$ &  $0.314$ & $0.558$ & $-0.570$\\
$y$ &  $0.442$ & $0.442$ & $-1.134$\\
$z$ & $0$ & $0.244$ & $-0.421$\\
$u$ & $0$ & $0$ & $ -0.422$\\
\hline
\end{tabular}
\end{table}
For the Ishigami function, the sum of all $S_i$ and first-order interaction terms $S_{ij}$ equals 1, and only $S_x$, $S_y$ and $S_{xz}$ are nonzero. Since $S_y=S_{Ty}$, $y$ is not involved in interaction terms, while $S_z$ does not have an effect by itself. The term $S_{xz}$ contributes to both $S_{Tx}$ and $S_{Tz}$.
From the numerical values for the Sobol indices, summarized in Table \ref{tab:sobol}, we conclude that $y$ is most important by itself, while $x$ is most important when interactions with other variables are included.

The proposed estimator $H_\alpha^*$ behaves similar to $S_i$. It is approximately $-0.42$ in case of independence (exact value depends on $N$, higher $N$ leads to a smaller value, see Figure \ref{fig:unifNvaryFMST}), and drops below this value if the variables are dependent. The values of $H_\alpha^*$ in Table \ref{tab:sobol} lead to the  same conclusion as those of $S_i$, of $y$ being the most important here, followed by $x$. The estimator $H_\alpha^*$ does not have the sum property, nor does it compute estimates for interaction effects. This may be possible by MSTs on higher-dimensional combinations of input and output variables, as mentioned in Section \ref{sec:5}.

In practice, it might be beneficial to compute all the Sobol indices when the input variables are independent and not large in number. However, our prime interest in this study is in general cases with dependent input variables given to us in the form of datasets. In these cases, Sobol indices are difficult to calculate, especially if one wants to take the dependencies into account. By contrast, $H_\alpha^*$ can be computed easily for dependent variables. For dependent variables, the expressions for the Sobol indices become more comprehensive and the interpretation becomes less straightforward, due to the possibility of negative values for the indices.

\subsection{Sediment transport}
In this testcase, we consider data of physical quantities related to sediment transport in the sea. The data consists of $N=9628$ measurements of three variables, and simulations of one intermediate variable and one output variable obtained at the North Sea coast near the town of Noordwijk (the Netherlands), between June 19, 1984 and October 5, 1987. The measured variables are root-mean-square wave height $H_{rms}$ (m), peak wave period $T_p$ (s), wave direction $\theta$ ($^\circ$) (relative to shore-normal). The output is the cross-shore sediment transport $S$ (m$^2$/s) computed with the Unibest-TC module \cite{3Rue07,2Wal12,2Wal00}. Also available is the long-shore sediment transport $S_y$ (m$^2$/s) obtained from brute force simulation without computation of the bed changes. Depending on the application, it can be treated as an input or an output variable.

We explore the dependencies between the input variables $(H_{rms}, T_p, \theta)$ and output variables $(S_y, S)$. However, before discussing the results, we want to point out a phenomenon frequently occurring in real-world datasets. Because of the accuracy of the measuring device, there are only 60 different values for $T_p$. This makes the data, although rank-transformed including the multiples, differing substantially from uniform. We correct for this in the rank-transform by random positioning within groups with the same value.
Similarly, for $H_{rms}$ and $\theta$ only 379 and 353 unique values, respectively, are present in the dataset. This leads to a few observed $(H_{rms},\theta)$ pairs occurring twice, so that the distance between these pairs is zero. Some MST algorithms treat this as a non-existing edge, resulting in larger values of the MST length. We circumvent this by adding a small noise to the normalized data matrix  ($\sigma=10^{-6}$), such that the datapoints are not exactly at the same location anymore, although the effect is too small to affect the ordering by the rank-transform.

The resulting normalized lengths are given in Table \ref{tab:sediment}.
\begin{table}[ht!]
\centering
\caption{Values of $H_\alpha^*$ for the sediment transport testcase.}
\label{tab:sediment}
\begin{tabular}{|r|cccc|}
\hline
& $T_p$ (s) & $\theta$ ($^\circ$) & $S_y$ (m$^2$/s) & $S$ (m$^2$/s)\\
\hline
$H_{rms}$ (m)& -0.579 & -0.467 & -0.687 & -0.709 \\
$T_p$ (s) & - & -0.490 & -0.529 & -0.529 \\
$\theta$ ($^\circ$) & - & - & -0.598 & -0.455 \\
$S_y$ (m$^2$/s) & - & - & - & -1.256 \\
\hline
\end{tabular}
\end{table}
It is found that the strongest dependency is between $S_y$ and $S$. This is not a surprise, since the long-shore and cross-shore sediment transport are expected to be dependent. The weakest dependencies are between $\theta$ and $S$ and between   $H_{rms}$ and $\theta$. The latter  indicates that the root-mean-square wave height and wave direction are only slightly correlated.

\section{Discussion}\label{sec:discussion}
In this study we have proposed a novel method for quantifying dependencies in multivariate datasets by means of the minimum spanning tree (MST). The length of the MST is directly related to the R\'{e}nyi entropy of the data, which in turn is a suitable quantity to assess dependency. Our approach to assess the relative strength of dependencies via the length of the MST can be used as an aid for sensitivity analysis, as discussed in Section \ref{sec:5} and demonstrated with numerical examples in Section \ref{sec:tests}.

Constructing the MST is computationally expensive for large datasets. To reduce the computational cost, algorithms that approximate the exact MST can be considered. We have compared three such approximation algorithms (Section \ref{sec:methods}), and found the multilevel FMST algorithm due to \cite{2Zho15} to be the most accurate.

The method proposed here can be considered as an alternative to the widely used variance-based sensitivity analysis methods. Because it is based on entropy rather than variance, it is easier to incorporate the dependencies.
Furthermore, it is well-suited for SA in cases where the input distributions are unknown and only a sample (dataset) of the inputs is available.

We have tested our approach on the Ishigami function as well as on a real-world testcase involving sediment transport in the North Sea.
In both cases we obtained rankings of dependencies consistent with prior knowledge or heuristic understanding.
Moreover, in the Ishigami test case we included an example with independent input variables, making it possible to compare our approach with the analysis based on Sobol indices.
The results from our proposed method were consistent with those from the Sobol indices.

Altogether, the approach proposed here is a suitable and useful method to quantify dependencies and to aid in the performance of SA.
It gives consistent results and remains computationally feasible even for large datasets by using the FMST algorithm.
A limitation of the method as presented in this paper is that it does not account for interaction effects, see the discussion in Section \ref{sec:5}.
However, as also briefly discussed in Section \ref{sec:5}, we anticipate that this method can be generalized by constructing MSTs in three (or more) dimensions, making it possible to account for interaction effects as well. We intend to explore this generalization in a future study.

\section*{Acknowledgements}
We would like to thank Bruna de Queiroz and Deltares for providing the sediment transport data. This research is part of the EUROS program, which is supported by NWO domain Applied and Engineering Sciences and partly funded by the Ministry of Economic Affairs.

\bibliographystyle{elsarticle-num}
\bibliography{references}
\end{document}